\documentclass[11pt]{article}%
\usepackage{amsfonts}
\usepackage{amsmath}
\usepackage{hyperref}
\usepackage{amssymb}
\usepackage{graphicx}%
\providecommand{\U}[1]{\protect\rule{.1in}{.1in}}
\begin{document}

\title{Change point analysis of an exponential model based on Phi-divergence test-statistics: simulated critical points case}
\author{Batsidis, A.$^{1}$, Mart\'{\i}n, N.$^{2}$\thanks{Corresponding author, E-mail:
\href{mailto: nirian.martin@uc3m.es}{nirian.martin@uc3m.es}.}, Pardo, L.$^{3}$
and Zografos, K.$^{1}$\\$^{1}${\small Dept. Mathematics, University of Ioannina, Greece}\\$^{2}${\small Dept. Statistics, Carlos III University of Madrid, Spain}\\$^{3}${\small Dept. Statistics and O.R., Complutense University of Madrid,
Spain}}
\maketitle

\begin{abstract}
Recently Batsidis \textit{et al.} (2011) have presented a new procedure based on divergence measures for testing the hypothesis of the existence of a change point in exponential populations. A simulation study was carried out, in this paper, using the asymptotic critical points obtained from the asymptotic distribution of the new test statistics introduced there. The main purpose of this paper is to study the behavior of the test statistics introduced in the cited paper of Batsidis \textit{et al.} (2011), using simulated critical points.

\end{abstract}


\medskip


\section{Introduction\label{sec0}}

In a previous paper Batsidis \textit{et al.} (2011) introduced new test statistics for
a change point in a sequence of independent exponentially distributed random
variables and studied their asymptotic distribution. Based on the asymptotic
critical point they presented a simulation study in order to analyze the
behavior of the new test statistics in relation with the Likelihood Ratio Test
(LRT) introduced and studied in Worsley (1986) and Haccou \textit{et al.}
(1985, 1988). The simulation study is based in obtaining the error of type I
as well as the power using the asymptotic critical points.

In this paper we shall study the behavior of the test statistics introduced in
the cited paper of Batsidis \textit{et al.} (2011) but using simulated critical points instead of asymptotic critical points.
In Section 2 we shall outline the different test statistics used for this purpose, while in Section 3 a
simulation study is carried out.

\section{Families of test statistics}

Let $X_{1},...,X_{K}$ be a sequence of $K$ independent exponential random
variables, with density function $f(x_{i},\theta_{i})=\theta_{i}exp\left(
-\theta_{i}x_{i}\right)  $, $i=1,...,K,$ respectively, where $\theta_{i}$ and
$x_{i}$ are positive real numbers. The change point problem is to test
hypothesis about the equality of means $1/\theta_{i}$, $i=1,...,K$, or what is
equivalent about the equality of the parameters $\theta_{i},$ $i=1,...,K,$
\begin{equation}
H_{0}:\text{ }\theta_{1}=\theta_{2}=...=\theta_{K}\text{ (}=\theta_{0}\text{,
}\theta_{0}\text{ unknown),}\label{A}%
\end{equation}
versus the alternative
\begin{equation}
H_{1}:\text{ }\theta_{1}=...=\theta_{k_{1}}\neq\theta_{k_{1}+1}=...=\theta
_{k_{2}}\neq...\neq\theta_{k_{q-1}+1}=...=\theta_{k_{q}}=\theta_{K},\label{B}%
\end{equation}
where $q$, $1\leq q\leq K$, is the unknown number of changes and $k_{1}%
,k_{2},...,k_{q}$ are the unknown positions of the change points.

Similar to the classical literature of change point analysis we just need to
test the single change point problem by means of the binary segmentation procedure.
This procedure was proposed by Vostrikova (1981). Based on this procedure we
just need to test the single change point hypothesis and then to repeat the
procedure for each subsequence. Therefore initially, we test for no change
point versus one change point, that is, we test the null hypothesis
\begin{equation}
H_{0}:X_{i}\text{ are described by }f_{\theta_{0}}(x)=\theta_{0}\exp
(-\theta_{0}x),\text{ }x>0\text{, }\theta_{0}>0\,;\text{ }i=1,...,K\text{ and
}\theta_{0}\text{ unknown,} \label{C}%
\end{equation}
versus the alternative
\begin{eqnarray}
H_{1}  &  :X_{i}\text{ are described by }f_{\theta_{0}}(x)=\theta_{0}%
\exp(-\theta_{0}x)\text{, }x>0\text{, }\theta_{0}>0\,;\text{ }i=1,...,k\text{
and }\label{D}\\
& X_{i}\text{ are described by }f_{\theta_{1}}(x)=\theta_{1}\exp(-\theta
_{1}x)\text{, }x>0\text{, }\theta_{1}>0\,;\text{ }i=k+1,...,K.\nonumber
\end{eqnarray}

Here, $k$ is the unknown location of the single change point. If $H_{0}$ is
not rejected, then the procedure is finished and there is no change point. If
$H_{0}$ is rejected, then there is a change point and we continue with the
step 2. In the second step we test separately the two subsequences before and
after the change point found in the first step. In the sequel, we repeat these
two steps until no further subsequences have change points. At the end of the
procedure, the collection of change point locations found by the previous
steps constitute the set of the change points. Therefore, we will concentrate
in the sequel on the case of a single change point in the random sequence.

It is well known that the maximum likelihood estimator (MLE), $\widehat{\theta
}_{0,k}^{(K)},$ of the parameter $\theta_{0}$ in the exponential model, which is
based on the random sample $X_{1},...,X_{k}$ from $f_{\theta_{0}}%
(x)=\theta_{0}\exp(-\theta_{0}x)$, is given by
\[
\widehat{\theta}_{0,k}^{(K)}=\frac{1}{\overline{X}_{k,0}}%
\]
with $\overline{X}_{k,0}=\frac{1}{k}%
{\textstyle\sum\nolimits_{i=1}^{k}}
X_{i}$ and the MLE of $\theta_{1}$ in the exponential model, which is based on
the random sample $X_{k+1},...,X_{K}$ from $f_{\theta_{1}}(x)=\theta_{1}%
\exp(-\theta_{1}x)$, is given by
\[
\widehat{\theta}_{1,k}^{(K)}=\frac{1}{\overline{X}_{K-k,1}}%
\]
with $\overline{X}_{K-k,1}=\frac{1}{K-k}%
{\textstyle\sum\nolimits_{i=k+1}^{K}}
X_{i}.$ We are using the subscript "$0$" to declare that we refer to the first
$k$ observations and the subscript "$1$" to declare that we refer to the last
$K-k$ last observations.

We denote,%
\begin{equation}
N(\epsilon)=\left\{  k/k\in\left\{  1,...,K-1\right\}  \text{
and }\frac{k}{K}\in\left[  \varepsilon,1-\varepsilon\right]  ,\text{ with
}\varepsilon>0\text{ and small enough}\right\}  .\label{D1}%
\end{equation}
Batsidis \textit{et al.} (2011b) considered the following family of test statistics
\begin{equation}
^{\epsilon}T_{\phi_{\lambda}}^{(K)}=\left\{
\begin{array}
[c]{lcc}%
\max_{k\in N(\epsilon)}\frac{2k(K-k)}{K\lambda\left(  \lambda+1\right)
}\left\{  \frac{\overline{X}_{k,0}^{-\lambda}\overline{X}_{K-k,1}^{\lambda+1}%
}{\left(  \lambda+1\right)  \overline{X}_{K-k,1}-\lambda\overline{X}_{k,0}%
}-1\right\}   &,& \left(  \lambda+1\right)  \overline{X}%
_{K-k,1}-\lambda\overline{X}_{k,0}>0,\text{ }\lambda\neq0,-1\\
\max_{k\in N(\epsilon)}\frac{2k(K-k)}{K}\left\{  \ln\frac{\overline{X}%
_{K-k,1}}{\overline{X}_{k,0}}+\frac{\overline{X}_{k,0}}{\overline{X}_{K-k,1}%
}-1\right\}   &,& \lambda=0\\
\max_{k\in N(\epsilon)}\frac{2k(K-k)}{K}\left\{  \ln\frac{\overline{X}_{k,0}%
}{\overline{X}_{K-k,1}}+\frac{\overline{X}_{K-k,1}}{\overline{X}_{k,0}%
}-1\right\}   &,& \lambda=-1
\end{array}
\right.  \label{I2}%
\end{equation}
for testing (\ref{C}) versus (\ref{D}). We have to note that when $\lambda
\neq0$, $\lambda\neq-1$ the condition $\left(  \lambda+1\right)  \overline
{X}_{K-k,1}-\lambda\overline{X}_{k,0}>0$ should be satisfied in order to ensure
the existence of the divergence. Taking into account that $\overline
{X}_{K-k,1}>0$ and $\overline{X}_{k,0}>0$ we will restrict to values of
$\lambda$ in the interval $(-1,0)$, so as the previous mentioned condition
to be always satisfied.

Worsley (1986) and Haccou \textit{et al.} (1985, 1988) used maximum
likelihood methods in order to test for a change in a sequence of independent
exponential family random variables, with particular emphasis on the
exponential distribution. It was proved that if $LRT^{(K)}$ is minus twice the
log likelihood ratio, that is
\begin{equation}
LRT^{(K)}=\underset{k\in\{1,...,K-1\}}{\max}-2\log\frac{%
{\textstyle\prod\limits_{h=1}^{K}}
f_{\widehat{\theta}_{0,K}^{(K)}}(X_{h})}{%
{\textstyle\prod\limits_{h=1}^{k}}
f_{\widehat{\theta}_{0,k}^{(K)}}(X_{h})%
{\textstyle\prod\limits_{h=k+1}^{K}}
f_{\widehat{\theta}_{1,k}^{(K)}}(X_{h})}, \label{J1}%
\end{equation}
with $f_{\widehat{\theta}_{0,k}^{(K)}}(X_{h})=\widehat{\theta}_{0,k}^{(K)}%
\exp(-\widehat{\theta}_{0,k}^{(K)}X_{h})$, and $f_{\widehat{\theta}%
_{1,k}^{(K)}}(X_{h})=\widehat{\theta}_{1,k}^{(K)}\exp(-\widehat{\theta}%
_{1,k}^{(K)}X_{h})$, $k=1,...,K-1$, then its explicit expression is given by
\begin{equation}
LRT^{(K)}=2\underset{k\in\{1,...,K-1\}}{\max}\left(  k\log\frac{\widehat
{\theta}_{0,k}^{(K)}}{\widehat{\theta}_{0,K}^{(K)}}+(K-k)\log\frac
{\widehat{\theta}_{1,k}^{(K)}}{\widehat{\theta}_{0,K}^{(K)}}\right)  .
\label{K1}%
\end{equation}

We have to note that $LRT^{(K)}$ can be written in terms of the
Kullback-Leibler divergence measure, as follows
\[
LRT^{(K)}=2K\underset{k\in\{1,...,K-1\}}{\max}\left(  \frac{k}{K}%
D_{Kull}\left(  f_{\widehat{\theta}_{0,k}^{(K)}},f_{\widehat{\theta}%
_{0,K}^{(K)}}\right)  +\frac{K-k}{K}D_{Kull}\left(  f_{\widehat{\theta}%
_{1,k}^{(K)}},f_{\widehat{\theta}_{0,K}^{(K)}}\right)  \right)  ,
\]
since in the special case of two exponential distributions with parameters
$\theta$ and $\theta^{\prime}$ it is easily obtained that
\[
D_{Kull}\left(  f_{\theta},f_{\theta^{\prime}}\right)  =\log\frac{\theta
}{\theta^{\prime}}+\frac{\theta^{\prime}}{\theta}-1.
\]

It is very interesting to observe that in many statistical problems LRT is
obtained directly from the Kullback-Leibler divergence but we can see that
in the case of the change point detection in the exponential model it is not possible. We can only express it as
a linear combination of two Kullback-Leibler divergences. This point is very
important because, based on it, the LRT is not a member of the
family of test statistics considered in (\ref{I2}). Batsidis \textit{et al.} (2011),
inspired on the third test-statistic proposed in Horv\'{a}th and Serbinowska
(1995, p. 373), which tends usually to behave much more better than the likelihood
ratio test statistic, considered a modified likelihood ratio test as a
weighted sum of two Kullback-Leibler divergences. The test statistic is given
by the following formula:
\begin{align}
S^{(K)} &  =2\underset{k\in\{1,...,K-1\}}{\max}\frac{k(K-k)}{K}\left(
\frac{k}{K}D_{Kull}\left(  f_{\widehat{\theta}_{0,k}^{(K)}},f_{\widehat
{\theta}_{0,K}^{(K)}}\right)  +\frac{K-k}{K}D_{Kull}\left(  f_{\widehat
{\theta}_{1,k}^{(K)}},f_{\widehat{\theta}_{0,K}^{(K)}}\right)  \right)
\nonumber\\
&  =2\underset{k\in\{1,...,K-1\}}{\max}\frac{k(K-k)}{K}\left(  \frac{k}{K}%
\log\frac{\widehat{\theta}_{0,k}^{(K)}}{\widehat{\theta}_{0,K}^{(K)}}%
+\frac{K-k}{K}\log\frac{\widehat{\theta}_{1,k}^{(K)}}{\widehat{\theta}%
_{0,K}^{(K)}}\right)  .\label{test1}%
\end{align}
In the next section we study the behavior of $^{\epsilon}T_{\phi_{\lambda}%
}^{(K)},$ $LRT^{(K)}$ and $S^{(K)}$ when simulated critical values are used, instead of the asymptotic critical points considered in Batsidis \textit{et al.} (2011).

\section{Simulation Study\label{sec3}}

The aim of this section is the evaluation of the simulated  critical values of
the different test statistics presented in Section 2, subject to the
assumption that there is no change point. In order to obtain the critical values for
the significance levels $0.1,0.05$ and $0.01$, we simulate $B=5000$ data sets of
sample size $K=40,50,60,64,100,200,300,400,500$ from the standard exponential
distribution. For each sample-data set, we calculate the test statistics
$^{0.05}T_{\phi_{\lambda}}^{(K)}$, for eleven values of $\lambda$, $\lambda
=-0.1(0.1)1$, as well as the test statistics $\widetilde
{LRT}^{(K)}=a(K)\sqrt{LRT^{(K)}}-b(K)$ and $S^{(K)}$. Therefore, 5000 values
of these statistics are obtained and the estimated $0.1$, $0.05$ and $0.01$
critical values are given in Table 1.

From Table 1 it can been seen that the estimated critical values for the
various significance levels are very different from the asymptotic critical
values (see Haccou \textit{et al.} (1985, p.10) for a similar conclusion). The
disadvantage of not using the limiting distribution is that specialized
extensive tables of simulated critical values for a simulation study when the
sample size does not exactly match the tabled values given are needed (see
Srivastava and Hui (1987)). However as for instance Romeu and Ozturk (1993)
pointed out using the limiting distribution implies loss of power.

\begin{table}[h]
\caption{Simulated critical values of the test statistics based on 5000
replications.}
\begin{center}
{\tiny
\begin{tabular}
[c]{llccccccccccccc}
&  & \multicolumn{11}{c}{$^{0.05}T_{\phi_{\lambda}}^{(K)}$} & $\widetilde
{LRT}^{(K)}$ & $S^{(K)}$\\
&  & \multicolumn{11}{c}{$\lambda$} &  & \\
$\alpha$ & $K$ & $-1$ & $-0.9$ & $-0.8$ & $-0.7$ & $-0.6$ & $-0.5$ & $-0.4$ &
$-0.3$ & $-0.2$ & $-0.1$ & $0$ &  & \\\hline
& $40$ & \multicolumn{1}{l}{11.5810} & \multicolumn{1}{l}{9.5676} &
\multicolumn{1}{l}{8.5411} & \multicolumn{1}{l}{7.9334} &
\multicolumn{1}{l}{7.6490} & \multicolumn{1}{l}{7.4920} &
\multicolumn{1}{l}{7.5585} & \multicolumn{1}{l}{7.8853} &
\multicolumn{1}{l}{8.4269} & \multicolumn{1}{l}{9.5196} &
\multicolumn{1}{l}{11.3937} & 2.3190 & 1.3025\\
& $50$ & \multicolumn{1}{l}{9.6585} & \multicolumn{1}{l}{8.5539} &
\multicolumn{1}{l}{7.9585} & \multicolumn{1}{l}{7.5286} &
\multicolumn{1}{l}{7.2882} & \multicolumn{1}{l}{7.2120} &
\multicolumn{1}{l}{7.3280} & \multicolumn{1}{l}{7.5125} &
\multicolumn{1}{l}{7.8517} & \multicolumn{1}{l}{8.5151} &
\multicolumn{1}{l}{9.5691} & 2.2360 & 1.3055\\
& $60$ & \multicolumn{1}{l}{9.8676} & \multicolumn{1}{l}{8.8737} &
\multicolumn{1}{l}{8.2602} & \multicolumn{1}{l}{7.8641} &
\multicolumn{1}{l}{7.6617} & \multicolumn{1}{l}{7.6264} &
\multicolumn{1}{l}{7.6725} & \multicolumn{1}{l}{7.8283} &
\multicolumn{1}{l}{8.1709} & \multicolumn{1}{l}{8.8069} &
\multicolumn{1}{l}{9.8090} & 2.2703 & 1.3030\\
& $64$ & \multicolumn{1}{l}{8.9705} & \multicolumn{1}{l}{8.2860} &
\multicolumn{1}{l}{7.8496} & \multicolumn{1}{l}{7.5334} &
\multicolumn{1}{l}{7.3362} & \multicolumn{1}{l}{7.3376} &
\multicolumn{1}{l}{7.3600} & \multicolumn{1}{l}{7.4888} &
\multicolumn{1}{l}{7.7334} & \multicolumn{1}{l}{8.1398} &
\multicolumn{1}{l}{8.7467} & 2.2961 & 1.2984\\
$0.1$ & $100$ & \multicolumn{1}{l}{8.7360} & \multicolumn{1}{l}{8.2728} &
\multicolumn{1}{l}{7.9100} & \multicolumn{1}{l}{7.7443} &
\multicolumn{1}{l}{7.6026} & \multicolumn{1}{l}{7.5428} &
\multicolumn{1}{l}{7.5825} & \multicolumn{1}{l}{7.7012} &
\multicolumn{1}{l}{7.9344} & \multicolumn{1}{l}{8.2489} &
\multicolumn{1}{l}{8.7130} & \multicolumn{1}{l}{2.2971} &
\multicolumn{1}{l}{1.3324}\\
& $200$ & \multicolumn{1}{l}{8.5279} & \multicolumn{1}{l}{8.3250} &
\multicolumn{1}{l}{8.2117} & \multicolumn{1}{l}{8.0575} &
\multicolumn{1}{l}{7.9950} & \multicolumn{1}{l}{7.9180} &
\multicolumn{1}{l}{7.9329} & \multicolumn{1}{l}{8.0188} &
\multicolumn{1}{l}{8.1164} & \multicolumn{1}{l}{8.2670} &
\multicolumn{1}{l}{8.6134} & \multicolumn{1}{l}{2.3153} &
\multicolumn{1}{l}{1.4177}\\
& $300$ & \multicolumn{1}{l}{8.3041} & \multicolumn{1}{l}{8.1771} &
\multicolumn{1}{l}{8.0462} & \multicolumn{1}{l}{7.9215} &
\multicolumn{1}{l}{7.9077} & \multicolumn{1}{l}{7.9109} &
\multicolumn{1}{l}{7.9861} & \multicolumn{1}{l}{7.9804} &
\multicolumn{1}{l}{8.0459} & \multicolumn{1}{l}{8.1457} &
\multicolumn{1}{l}{8.3344} & 2.3649 & 1.3813\\
& $400$ & \multicolumn{1}{l}{8.2505} & \multicolumn{1}{l}{8.1697} &
\multicolumn{1}{l}{8.0972} & \multicolumn{1}{l}{8.0184} &
\multicolumn{1}{l}{7.9514} & \multicolumn{1}{l}{7.9289} &
\multicolumn{1}{l}{7.9042} & \multicolumn{1}{l}{7.9476} &
\multicolumn{1}{l}{7.9841} & \multicolumn{1}{l}{8.0687} &
\multicolumn{1}{l}{82334} & 2.3606 & 1.4142\\
& $500$ & \multicolumn{1}{l}{8.2393} & \multicolumn{1}{l}{8.1955} &
\multicolumn{1}{l}{8.1057} & \multicolumn{1}{l}{8.0541} &
\multicolumn{1}{l}{8.0201} & \multicolumn{1}{l}{7.9773} &
\multicolumn{1}{l}{7.9877} & \multicolumn{1}{l}{7.9788} &
\multicolumn{1}{l}{8.0268} & \multicolumn{1}{l}{8.1105} &
\multicolumn{1}{l}{8.2279} & \multicolumn{1}{l}{2.3646} &
\multicolumn{1}{l}{1.4229}\\
& $\infty$ & \multicolumn{1}{l}{8.31} & \multicolumn{1}{l}{8.31} &
\multicolumn{1}{l}{8.31} & \multicolumn{1}{l}{8.31} & \multicolumn{1}{l}{8.31}
& \multicolumn{1}{l}{8.31} & \multicolumn{1}{l}{8.31} &
\multicolumn{1}{l}{8.31} & \multicolumn{1}{l}{8.31} & \multicolumn{1}{l}{8.31}
& \multicolumn{1}{l}{8.31} & \multicolumn{1}{l}{2.9435} &
\multicolumn{1}{l}{1.4978}\\\hline
& $40$ & \multicolumn{1}{l}{17.6334} & \multicolumn{1}{l}{13.1196} &
\multicolumn{1}{l}{11.0430} & \multicolumn{1}{l}{9.8748} &
\multicolumn{1}{l}{9.3129} & \multicolumn{1}{l}{9.1339} &
\multicolumn{1}{l}{9.1892} & \multicolumn{1}{l}{9.7084} &
\multicolumn{1}{l}{10.7864} & \multicolumn{1}{l}{12.7154} &
\multicolumn{1}{l}{17.6140} & 2.7241 & 1.6569\\
& $50$ & \multicolumn{1}{l}{13.5357} & \multicolumn{1}{l}{11.3328} &
\multicolumn{1}{l}{9.9970} & \multicolumn{1}{l}{9.2835} &
\multicolumn{1}{l}{8.9350} & \multicolumn{1}{l}{8.8046} &
\multicolumn{1}{l}{8.8077} & \multicolumn{1}{l}{9.1188} &
\multicolumn{1}{l}{9.7586} & \multicolumn{1}{l}{11.0303} &
\multicolumn{1}{l}{13.1107} & 2.6826 & 1.6648\\
& $60$ & \multicolumn{1}{l}{13.9701} & \multicolumn{1}{l}{11.5162} &
\multicolumn{1}{l}{10.3217} & \multicolumn{1}{l}{9.6691} &
\multicolumn{1}{l}{9.2967} & \multicolumn{1}{l}{9.2192} &
\multicolumn{1}{l}{9.3095} & \multicolumn{1}{l}{9.7703} &
\multicolumn{1}{l}{10.5004} & \multicolumn{1}{l}{11.7857} &
\multicolumn{1}{l}{13.9239} & 2.6850 & 1.6546\\
& $64$ & \multicolumn{1}{l}{12.2527} & \multicolumn{1}{l}{10.7728} &
\multicolumn{1}{l}{9.9618} & \multicolumn{1}{l}{9.4196} &
\multicolumn{1}{l}{9.0943} & \multicolumn{1}{l}{9.0235} &
\multicolumn{1}{l}{9.0149} & \multicolumn{1}{l}{9.2659} &
\multicolumn{1}{l}{9.8645} & \multicolumn{1}{l}{10.6607} &
\multicolumn{1}{l}{12.1021} & \multicolumn{1}{l}{2.7308} &
\multicolumn{1}{l}{1.6309}\\
$0.05$ & $100$ & \multicolumn{1}{l}{11.4735} & \multicolumn{1}{l}{10.5057} &
\multicolumn{1}{l}{9.8783} & \multicolumn{1}{l}{9.3770} &
\multicolumn{1}{l}{9.1050} & \multicolumn{1}{l}{9.0232} &
\multicolumn{1}{l}{9.1798} & \multicolumn{1}{l}{9.4268} &
\multicolumn{1}{l}{9.7724} & \multicolumn{1}{l}{10.4078} &
\multicolumn{1}{l}{11.1852} & \multicolumn{1}{l}{2.7164} &
\multicolumn{1}{l}{1.6676}\\
& $200$ & \multicolumn{1}{l}{10.6694} & \multicolumn{1}{l}{10.1903} &
\multicolumn{1}{l}{9.9174} & \multicolumn{1}{l}{9.6909} &
\multicolumn{1}{l}{9.6413} & \multicolumn{1}{l}{9.5670} &
\multicolumn{1}{l}{9.7049} & \multicolumn{1}{l}{9.8352} &
\multicolumn{1}{l}{10.0539} & \multicolumn{1}{l}{10.3680} &
\multicolumn{1}{l}{10.8759} & \multicolumn{1}{l}{2.8402} &
\multicolumn{1}{l}{1.7269}\\
& $300$ & \multicolumn{1}{l}{10.3158} & \multicolumn{1}{l}{10.0808} &
\multicolumn{1}{l}{9.8837} & \multicolumn{1}{l}{9.7218} &
\multicolumn{1}{l}{9.5948} & \multicolumn{1}{l}{9.5428} &
\multicolumn{1}{l}{9.5377} & \multicolumn{1}{l}{9.6841} &
\multicolumn{1}{l}{9.8029} & \multicolumn{1}{l}{9.9631} &
\multicolumn{1}{l}{10.2826} & 2.8237 & 1.7393\\
& $400$ & \multicolumn{1}{l}{10.1254} & \multicolumn{1}{l}{9.9236} &
\multicolumn{1}{l}{9.8319} & \multicolumn{1}{l}{9.6970} &
\multicolumn{1}{l}{9.5952} & \multicolumn{1}{l}{9.5807} &
\multicolumn{1}{l}{9.5810} & \multicolumn{1}{l}{9.5850} &
\multicolumn{1}{l}{9.7662} & \multicolumn{1}{l}{9.9785} &
\multicolumn{1}{l}{10.1668} & 2.8461 & 1.7573\\
& $500$ & \multicolumn{1}{l}{10.1375} & \multicolumn{1}{l}{9.9245} &
\multicolumn{1}{l}{9.7938} & \multicolumn{1}{l}{9.7101} &
\multicolumn{1}{l}{9.6954} & \multicolumn{1}{l}{9.6010} &
\multicolumn{1}{l}{9.5865} & \multicolumn{1}{l}{9.6099} &
\multicolumn{1}{l}{9.6936} & \multicolumn{1}{l}{9.8310} &
\multicolumn{1}{l}{10.0460} & \multicolumn{1}{l}{2.8415} &
\multicolumn{1}{l}{1.7784}\\
& $\infty$ & \multicolumn{1}{l}{9.90} & \multicolumn{1}{l}{9.90} &
\multicolumn{1}{l}{9.90} & \multicolumn{1}{l}{9.90} & \multicolumn{1}{l}{9.90}
& \multicolumn{1}{l}{9.90} & \multicolumn{1}{l}{9.90} &
\multicolumn{1}{l}{9.90} & \multicolumn{1}{l}{9.90} & \multicolumn{1}{l}{9.90}
& \multicolumn{1}{l}{9.90} & \multicolumn{1}{l}{3.6633} &
\multicolumn{1}{l}{1.8444}\\\hline
& $40$ & \multicolumn{1}{l}{50.5550} & \multicolumn{1}{l}{24.4875} &
\multicolumn{1}{l}{17.3032} & \multicolumn{1}{l}{14.1739} &
\multicolumn{1}{l}{12.7427} & \multicolumn{1}{l}{12.1567} &
\multicolumn{1}{l}{12.5086} & \multicolumn{1}{l}{13.6526} &
\multicolumn{1}{l}{16.3834} & \multicolumn{1}{l}{22.9197} &
\multicolumn{1}{l}{44.9917} & 3.5780 & 2.4148\\
& $50$ & \multicolumn{1}{l}{27.0400} & \multicolumn{1}{l}{19.4354} &
\multicolumn{1}{l}{15.5167} & \multicolumn{1}{l}{13.6776} &
\multicolumn{1}{l}{12.6433} & \multicolumn{1}{l}{12.3613} &
\multicolumn{1}{l}{12.4914} & \multicolumn{1}{l}{13.2924} &
\multicolumn{1}{l}{15.4007} & \multicolumn{1}{l}{19.0322} &
\multicolumn{1}{l}{28.8758} & 3.5656 & 2.4720\\
& $60$ & \multicolumn{1}{l}{29.6808} & \multicolumn{1}{l}{19.8103} &
\multicolumn{1}{l}{15.9025} & \multicolumn{1}{l}{13.7726} &
\multicolumn{1}{l}{12.8346} & \multicolumn{1}{l}{12.6437} &
\multicolumn{1}{l}{12.7626} & \multicolumn{1}{l}{13.8295} &
\multicolumn{1}{l}{15.7556} & \multicolumn{1}{l}{19.7567} &
\multicolumn{1}{l}{28.2007} & 3.5103 & 2.4652\\
& $64$ & \multicolumn{1}{l}{21.1749} & \multicolumn{1}{l}{16.8784} &
\multicolumn{1}{l}{14.3429} & \multicolumn{1}{l}{13.1193} &
\multicolumn{1}{l}{12.5352} & \multicolumn{1}{l}{12.5130} &
\multicolumn{1}{l}{12.7935} & \multicolumn{1}{l}{13.5557} &
\multicolumn{1}{l}{15.1971} & \multicolumn{1}{l}{17.5568} &
\multicolumn{1}{l}{22.5916} & 3.6179 & 2.4349\\
$0.01$ & $100$ & \multicolumn{1}{l}{21.0281} & \multicolumn{1}{l}{17.3015} &
\multicolumn{1}{l}{14.9784} & \multicolumn{1}{l}{13.6942} &
\multicolumn{1}{l}{12.8163} & \multicolumn{1}{l}{12.7267} &
\multicolumn{1}{l}{12.8364} & \multicolumn{1}{l}{13.3041} &
\multicolumn{1}{l}{14.1556} & \multicolumn{1}{l}{16.1152} &
\multicolumn{1}{l}{19.3112} & 3.4667 & 2.3155\\
& $200$ & \multicolumn{1}{l}{16.3509} & \multicolumn{1}{l}{15.0467} &
\multicolumn{1}{l}{14.1423} & \multicolumn{1}{l}{13.9136} &
\multicolumn{1}{l}{13.7878} & \multicolumn{1}{l}{13.4483} &
\multicolumn{1}{l}{13.6832} & \multicolumn{1}{l}{14.3627} &
\multicolumn{1}{l}{14.8642} & \multicolumn{1}{l}{15.8700} &
\multicolumn{1}{l}{17.5407} & \multicolumn{1}{l}{3.7306} &
\multicolumn{1}{l}{2.5576}\\
& $300$ & \multicolumn{1}{l}{15.8743} & \multicolumn{1}{l}{15.2039} &
\multicolumn{1}{l}{14.5089} & \multicolumn{1}{l}{13.9226} &
\multicolumn{1}{l}{13.9182} & \multicolumn{1}{l}{13.7411} &
\multicolumn{1}{l}{13.5544} & \multicolumn{1}{l}{13.5335} &
\multicolumn{1}{l}{13.9115} & \multicolumn{1}{l}{14.7973} &
\multicolumn{1}{l}{15.6554} & 3.7907 & 2.5829\\
& $400$ & \multicolumn{1}{l}{14.3163} & \multicolumn{1}{l}{14.1375} &
\multicolumn{1}{l}{13.7048} & \multicolumn{1}{l}{13.3444} &
\multicolumn{1}{l}{13.1254} & \multicolumn{1}{l}{13.2488} &
\multicolumn{1}{l}{13.4639} & \multicolumn{1}{l}{13.7315} &
\multicolumn{1}{l}{14.1345} & \multicolumn{1}{l}{14.6297} &
\multicolumn{1}{l}{15.1618} & 3.8562 & 2.5797\\
& $500$ & \multicolumn{1}{l}{14.8829} & \multicolumn{1}{l}{14.4170} &
\multicolumn{1}{l}{14.0629} & \multicolumn{1}{l}{13.6848} &
\multicolumn{1}{l}{13.3272} & \multicolumn{1}{l}{13.2033} &
\multicolumn{1}{l}{13.2050} & \multicolumn{1}{l}{13.2802} &
\multicolumn{1}{l}{13.6211} & \multicolumn{1}{l}{13.9468} &
\multicolumn{1}{l}{14.3068} & \multicolumn{1}{l}{3.8130} &
\multicolumn{1}{l}{2.5964}\\
& $\infty$ & \multicolumn{1}{l}{13.45} & \multicolumn{1}{l}{13.45} &
\multicolumn{1}{l}{13.45} & \multicolumn{1}{l}{13.45} &
\multicolumn{1}{l}{13.45} & \multicolumn{1}{l}{13.45} &
\multicolumn{1}{l}{13.45} & \multicolumn{1}{l}{13.45} &
\multicolumn{1}{l}{13.45} & \multicolumn{1}{l}{13.45} &
\multicolumn{1}{l}{13.45} & \multicolumn{1}{l}{5.2933} &
\multicolumn{1}{l}{2.6491}%
\end{tabular}
}
\end{center}
\end{table}

The type I error rate is an essential characteristic of the performance of a
test statistic. In the sequel we present the results of a Monte Carlo study on
the type I error rates, by considering data sampled from a standard
exponential distribution with no change point. In the study 5000 data sets
with different sample size $K=40,50,60,64,100,200,300,500$ were generated. The
simulated results on the type I error rates of $^{0.05}T_{\phi_{\lambda}%
}^{(K)}$, $\lambda
=-0.1(0.1)1$, as well as the test
statistics $\widetilde{LRT}^{(K)}$ and $S^{(K)}$ for testing the existence of
a change point, for significance level $\alpha=0.1$, $0.05$, $0.01$, are
presented in Table 2, when the simulated critical values of Table 1 are used
respectively. In this framework the test statistics were calculated for all
5000 data sets and compared to the appropriate simulated critical value given
in Table 1. The type I error rate is estimated as the proportion of rejections
of the hypothesis of no change point in each situation. So in Table 2 the
empirical sizes, the proportion of times that the null hypothesis is rejected
when all data are distributed according to the standard exponential
distribution are given.

Since these type I error rates were estimated using Monte Carlo simulations,
they are not free of error. So in order to decide if a test is accurate or
not, in a similar manner with that of Cardoso de Oliveira and Ferreira (2010), we
apply the exact binomial test for the null hypothesis $H_{0}: \alpha=0.1
(0.05, 0.01$, respectively) versus the alternative $H_{1}:\alpha\neq0.1 (0.05,
0.01$, respectively), with a significance level of $0.01$. In this frame, if
the simulated (observed) type I error rates were not significantly different
from $\alpha$ the test would be considered accurate. Otherwise, the test would
be considered conservative or liberal. In Table 2 we indicate with $*$ the
significantly different simulated type I error rates from the nominal level
$\alpha$. Based on Table 2 we conclude that in almost all cases presented
here the test statistics considered are accurate when the simulated critical
values of the Table 1 were used.

\begin{table}[h]
\caption{Empirical sizes based on 5000 replications with simulated critical
values.}
\begin{center}
{\tiny
\begin{tabular}
[c]{llccccccccccccc}
&  & \multicolumn{11}{c}{$^{0.05}T_{\phi_{\lambda}}^{(K)}$} & $\widetilde
{LRT}^{(K)}$ & $S^{(K)}$\\
&  & \multicolumn{11}{c}{$\lambda$} &  & \\
$\alpha$ & $K$ & $-1$ & $-0.9$ & $-0.8$ & $-0.7$ & $-0.6$ & $-0.5$ & $-0.4$ &
$-0.3$ & $-0.2$ & $-0.1$ & $0$ &  & \\\hline
& $40$ & \multicolumn{1}{l}{0.0968} & \multicolumn{1}{l}{0.0928} &
\multicolumn{1}{l}{0.0896} & \multicolumn{1}{l}{0.0932} &
\multicolumn{1}{l}{0.0952} & \multicolumn{1}{l}{0.0980} &
\multicolumn{1}{l}{0.0962} & \multicolumn{1}{l}{0.0960} &
\multicolumn{1}{l}{0.0930} & \multicolumn{1}{l}{0.0920} &
\multicolumn{1}{l}{0.0960} & \multicolumn{1}{l}{0.0928} &
\multicolumn{1}{l}{0.1022}\\
& $50$ & \multicolumn{1}{l}{0.0986} & \multicolumn{1}{l}{0.1042} &
\multicolumn{1}{l}{0.1018} & \multicolumn{1}{l}{0.1032} &
\multicolumn{1}{l}{0.1040} & \multicolumn{1}{l}{0.1054} &
\multicolumn{1}{l}{0.1026} & \multicolumn{1}{l}{0.1040} &
\multicolumn{1}{l}{0.1072} & \multicolumn{1}{l}{0.1038} &
\multicolumn{1}{l}{0.1062} & \multicolumn{1}{l}{0.1106} &
\multicolumn{1}{l}{0.1050}\\
& $60$ & \multicolumn{1}{l}{0.0968} & \multicolumn{1}{l}{0.0964} &
\multicolumn{1}{l}{0.0940} & \multicolumn{1}{l}{0.0940} &
\multicolumn{1}{l}{0.0970} & \multicolumn{1}{l}{0.0942} &
\multicolumn{1}{l}{0.0990} & \multicolumn{1}{l}{0.1008} &
\multicolumn{1}{l}{0.1032} & \multicolumn{1}{l}{0.0996} &
\multicolumn{1}{l}{0.1004} & \multicolumn{1}{l}{0.0992} &
\multicolumn{1}{l}{0.1018}\\
& $64$ & \multicolumn{1}{l}{0.1038} & \multicolumn{1}{l}{0.1006} &
\multicolumn{1}{l}{0.1010} & \multicolumn{1}{l}{0.1032} &
\multicolumn{1}{l}{0.1034} & \multicolumn{1}{l}{0.0986} &
\multicolumn{1}{l}{0.1004} & \multicolumn{1}{l}{0.1012} &
\multicolumn{1}{l}{0.1012} & \multicolumn{1}{l}{0.1026} &
\multicolumn{1}{l}{0.1026} & \multicolumn{1}{l}{0.0978} &
\multicolumn{1}{l}{0.1056}\\
$0.1$ & $100$ & \multicolumn{1}{l}{0.1048} & \multicolumn{1}{l}{0.1034} &
\multicolumn{1}{l}{0.1028} & \multicolumn{1}{l}{0.0976} &
\multicolumn{1}{l}{0.0998} & \multicolumn{1}{l}{0.0992} &
\multicolumn{1}{l}{0.0988} & \multicolumn{1}{l}{0.1004} &
\multicolumn{1}{l}{0.0970} & \multicolumn{1}{l}{0.0980} &
\multicolumn{1}{l}{0.0990} & 0.0960 & 0.1030\\
& $200$ & \multicolumn{1}{l}{0.0998} & \multicolumn{1}{l}{0.0970} &
\multicolumn{1}{l}{0.0946} & \multicolumn{1}{l}{0.0958} &
\multicolumn{1}{l}{0.0970} & \multicolumn{1}{l}{0.1014} &
\multicolumn{1}{l}{0.1030} & \multicolumn{1}{l}{0.1002} &
\multicolumn{1}{l}{0.1022} & \multicolumn{1}{l}{0.1022} &
\multicolumn{1}{l}{0.0982} & 0.1066 & 0.0958\\
& $300$ & \multicolumn{1}{l}{0.1000} & \multicolumn{1}{l}{0.1000} &
\multicolumn{1}{l}{0.1000} & \multicolumn{1}{l}{0.1000} &
\multicolumn{1}{l}{0.1000} & \multicolumn{1}{l}{0.1000} &
\multicolumn{1}{l}{0.1000} & \multicolumn{1}{l}{0.1000} &
\multicolumn{1}{l}{0.1000} & \multicolumn{1}{l}{0.1000} &
\multicolumn{1}{l}{0.1000} & \multicolumn{1}{l}{0.1000} &
\multicolumn{1}{l}{0.1000}\\
& $500$ & \multicolumn{1}{l}{0.1002} & \multicolumn{1}{l}{0.0972} &
\multicolumn{1}{l}{0.0982} & \multicolumn{1}{l}{0.0988} &
\multicolumn{1}{l}{0.0988} & \multicolumn{1}{l}{0.0986} &
\multicolumn{1}{l}{0.0986} & \multicolumn{1}{l}{0.1012} &
\multicolumn{1}{l}{0.1006} & \multicolumn{1}{l}{0.0990} &
\multicolumn{1}{l}{0.1000} & \multicolumn{1}{l}{0.1016} &
\multicolumn{1}{l}{0.1030}\\\hline
& $40$ & 0.0458 & 0.0472 & 0.0454 & 0.0424 & 0.0430 & 0.0450 & 0.0480 &
0.0484 & 0.0494 & 0.0490 & 0.0426 & 0.0460 & 0.0456\\
& $50$ & 0.0478 & 0.0494 & 0.0506 & 0.0524 & 0.0528 & 0.0530 & 0.0556 &
0.0568 & 0.0580 & 0.0546 & 0.0544 & 0.0516 & 0.0486\\
& $60$ & 0.0468 & 0.0490 & 0.0490 & 0.0480 & 0.0472 & 0.0482 & 0.0512 &
0.0478 & 0.0472 & 0.0454 & 0.0494 & 0.0536 & 0.0484\\
& $64$ & \multicolumn{1}{l}{0.0470} & \multicolumn{1}{l}{0.0476} &
\multicolumn{1}{l}{0.0478} & \multicolumn{1}{l}{0.0484} &
\multicolumn{1}{l}{0.0482} & \multicolumn{1}{l}{0.0506} &
\multicolumn{1}{l}{0.0540} & \multicolumn{1}{l}{0.0536} &
\multicolumn{1}{l}{0.0522} & \multicolumn{1}{l}{0.0536} &
\multicolumn{1}{l}{0.0518} & 0.0470 & 0.0514\\
$0.05$ & $100$ & \multicolumn{1}{l}{0.0508} & \multicolumn{1}{l}{0.0496} &
\multicolumn{1}{l}{0.0512} & \multicolumn{1}{l}{0.0542} &
\multicolumn{1}{l}{0.0546} & \multicolumn{1}{l}{0.0526} &
\multicolumn{1}{l}{0.0500} & \multicolumn{1}{l}{0.0492} &
\multicolumn{1}{l}{0.0516} & \multicolumn{1}{l}{0.0510} &
\multicolumn{1}{l}{0.0532} & 0.0512 & 0.0448\\
& $200$ & \multicolumn{1}{l}{0.0532} & \multicolumn{1}{l}{0.0532} &
\multicolumn{1}{l}{0.0540} & \multicolumn{1}{l}{0.0534} &
\multicolumn{1}{l}{0.0498} & \multicolumn{1}{l}{0.0506} &
\multicolumn{1}{l}{0.0474} & \multicolumn{1}{l}{0.0466} &
\multicolumn{1}{l}{0.0458} & \multicolumn{1}{l}{0.0458} &
\multicolumn{1}{l}{0.0442} & 0.0484 & 0.0532\\
& $300$ & 0.0500 & 0.0500 & 0.0500 & 0.0500 & 0.0500 & 0.0500 & 0.0500 &
0.0500 & 0.0500 & 0.0500 & 0.0500 & 0.0502 & 0.0500\\
& $500$ & \multicolumn{1}{l}{0.0476} & \multicolumn{1}{l}{0.0482} &
\multicolumn{1}{l}{0.0496} & \multicolumn{1}{l}{0.0504} &
\multicolumn{1}{l}{0.0496} & \multicolumn{1}{l}{0.0504} &
\multicolumn{1}{l}{0.0502} & \multicolumn{1}{l}{0.0488} &
\multicolumn{1}{l}{0.0490} & \multicolumn{1}{l}{0.0486} &
\multicolumn{1}{l}{0.0474} & 0.0488 & 0.0526\\\hline
& $40$ & 0.0076 & 0.0072 & 0.0076 & 0.0088 & 0.0092 & 0.0098 & 0.0092 &
0.0070 & 0.0070 & 0.0088 & 0.0094 & 0.0104 & 0.0116\\
& $50$ & 0.0114 & 0.0120 & 0.0116 & 0.0114 & 0.0114 & 0.0118 & 0.0134 &
0.0144$^{*}$ & 0.0128 & 0.0136 & 0.0106 & 0.0116 & 0.0076\\
& $60$ & 0.0092 & 0.0106 & 0.0100 & 0.0110 & 0.0108 & 0.0090 & 0.0092 &
0.0096 & 0.0098 & 0.0100 & 0.0114 & 0.0100 & 0.0098\\
& $64$ & \multicolumn{1}{l}{0.0104} & \multicolumn{1}{l}{0.0096} &
\multicolumn{1}{l}{0.0110} & \multicolumn{1}{l}{0.0114} &
\multicolumn{1}{l}{0.0120} & \multicolumn{1}{l}{0.0116} &
\multicolumn{1}{l}{0.0120} & \multicolumn{1}{l}{0.0140$^{*}$} &
\multicolumn{1}{l}{0.0134} & \multicolumn{1}{l}{0.0136} &
\multicolumn{1}{l}{0.0154$^{*}$} & 0.0080 & 0.0102\\
$0.01$ & $100$ & \multicolumn{1}{l}{0.0102} & \multicolumn{1}{l}{0.0104} &
\multicolumn{1}{l}{0.0042$^{*}$} & \multicolumn{1}{l}{0.0110} &
\multicolumn{1}{l}{0.0124} & \multicolumn{1}{l}{0.0120} &
\multicolumn{1}{l}{0.0118} & \multicolumn{1}{l}{0.0128} &
\multicolumn{1}{l}{0.0154$^{*}$} & \multicolumn{1}{l}{0.0144$^{*}$} &
\multicolumn{1}{l}{0.0126} & 0.0140$^{*}$ & 0.0132\\
& $200$ & \multicolumn{1}{l}{0.0142$^{*}$} & \multicolumn{1}{l}{0.0144$^{*}$}
& \multicolumn{1}{l}{0.0138$^{*}$} & \multicolumn{1}{l}{0.0110} &
\multicolumn{1}{l}{0.0078} & \multicolumn{1}{l}{0.0082} &
\multicolumn{1}{l}{0.0078} & \multicolumn{1}{l}{0.0064$^{*}$} &
\multicolumn{1}{l}{0.0070} & \multicolumn{1}{l}{0.0074} &
\multicolumn{1}{l}{0.0072} & 0.0082 & 0.0078\\
& $300$ & 0.0100 & 0.0100 & 0.0100 & 0.0100 & 0.0100 & 0.0100 & 0.0100 &
0.0100 & 0.0100 & 0.0100 & 0.0100 & 0.0100 & 0.0100\\
& $500$ & \multicolumn{1}{l}{0.0096} & \multicolumn{1}{l}{0.0092} &
\multicolumn{1}{l}{0.0096} & \multicolumn{1}{l}{0.0092} &
\multicolumn{1}{l}{0.0102} & \multicolumn{1}{l}{0.0098} &
\multicolumn{1}{l}{0.0096} & \multicolumn{1}{l}{0.0094} &
\multicolumn{1}{l}{0.0088} & \multicolumn{1}{l}{0.0090} &
\multicolumn{1}{l}{0.0096} & 0.0094 & 0.0094
\end{tabular}
$^{*}$ indicates that the simulated level is significantly different from
$\alpha$. }
\end{center}
\end{table}

Our aim in the sequel is to estimate the power of the test statistics for
detecting the existence of a change point by means of Monte Carlo methods. In
this context, motivated by Haccou \textit{et al.} (1985, 1988), we apply the
test statistics to $B=5000$ samples, with different significance levels
$\alpha$, as well as with several combinations of $K$ and $k$, with $K$ the
sample size and $k$ the point where change is taken place. Point $k$ is
selected such as $[k]=\tau\times K$, where $0<\tau<1$ and $[\cdot]$ denotes
the integer part of a real number.

For each combination of sample size $K$ and $k$, 5000 data sets were generated
such as $X_{1},...X_{k}$ to be distributed according to an exponential
distribution with $\theta_{0}=1$ and $X_{k+1},...,X_{K}$ to be described by
other exponential distribution with parameter $\theta_{1}\neq1$. In this
framework we consider the following scenarios $\theta_{1}\in
\{5,4,3,2,1/2,1/3,1/4,1/5\}$ and we denote by $\rho$ the ratio $\frac
{\theta_{1}}{\theta_{0}}$. For each sample we calculate the value of the test
statistics and decide to reject or not the (false) null hypothesis, with
significance level $\alpha$, considering the appropriate simulated or
asymptotical respectively critical value. The power is obtained calculating
the proportion in times in 5000 Monte Carlo simulations that the (false) null
hypothesis is rejected considering the specified significance level based on
the simulated critical values of Table 1.

In the sequel, we will present in Tables 3-8 the results for $\alpha=0.05,
0.01$, $K=40,50,100,200$ and $\tau=0.2, 0.3,0.5$.

\begin{table}[ptb]
\caption{Empirical powers based on 5000 replications with simulated critical
values, when $\alpha=0.05$ and $\tau=0.2$.}
\begin{center}
{\tiny
\begin{tabular}
[c]{llccccccccccccc}
&  & \multicolumn{11}{c}{$^{0.05}T_{\phi_{\lambda}}^{(K)}$} & $\widetilde
{LRT}^{(K)}$ & $S^{(K)}$\\
&  & \multicolumn{11}{c}{$\lambda$} &  & \\
$K$ & \multicolumn{1}{r}{$\theta_{1}$} & $-1$ & $-0.9$ & $-0.8$ & $-0.7$ &
$-0.6$ & $-0.5$ & $-0.4$ & $-0.3$ & $-0.2$ & $-0.1$ & $0$ &  & \\\hline
$40$ & \multicolumn{1}{r}{$5$} & 0.0746 & 0.3402 & 0.5550 & 0.6984 & 0.7804 &
0.8250 & 0.8528 & 0.8668 & 0.8700 & 0.8648 & 0.8246 & 0.9132 & 0.8994\\
& \multicolumn{1}{r}{$4$} & 0.0300 & 0.1710 & 0.3492 & 0.5010 & 0.5944 &
0.6636 & 0.7072 & 0.7330 & 0.7362 & 0.7270 & 0.6630 & 0.8176 & 0.7850\\
& \multicolumn{1}{r}{$3$} & 0.0120 & 0.0574 & 0.1434 & 0.2456 & 0.3312 &
0.3946 & 0.4492 & 0.4808 & 0.4888 & 0.4756 & 0.4018 & 0.5922 & 0.5610\\
& \multicolumn{1}{r}{$2$} & 0.0112 & 0.0142 & 0.0286 & 0.0564 & 0.0918 &
0.1256 & 0.1562 & 0.1744 & 0.1798 & 0.1758 & 0.1440 & 0.2340 & 0.2228\\
& \multicolumn{1}{r}{$1/2$} & 0.2396 & 0.2568 & 0.2644 & 0.2598 & 0.2376 &
0.2060 & 0.1790 & 0.1310 & 0.0848 & 0.0534 & 0.0376 & 0.1786 & 0.1844\\
& \multicolumn{1}{r}{$1/3$} & 0.5322 & 0.5714 & 0.5776 & 0.5736 & 0.5466 &
0.5066 & 0.4498 & 0.3606 & 0.2514 & 0.1368 & 0.0452 & 0.4334 & 0.4122\\
& \multicolumn{1}{r}{$1/4$} & 0.7536 & 0.7896 & 0.8004 & 0.7962 & 0.7780 &
0.7394 & 0.6896 & 0.6062 & 0.4606 & 0.2758 & 0.0716 & 0.6796 & 0.6432\\
& \multicolumn{1}{r}{$1/5$} & 0.8886 & 0.9122 & 0.9190 & 0.9170 & 0.9050 &
0.8854 & 0.8546 & 0.7910 & 0.6714 & 0.4552 & 0.1436 & 0.8424 & 0.8116\\\hline
$50$ & \multicolumn{1}{r}{$5$} & 0.5046 & 0.7160 & 0.8222 & 0.8802 & 0.9100 &
0.9302 & 0.9452 & 0.9532 & 0.9570 & 0.9556 & 0.9512 & 0.9654 & 0.9558\\
& \multicolumn{1}{r}{$4$} & 0.2654 & 0.4786 & 0.6310 & 0.7206 & 0.7780 &
0.8154 & 0.8454 & 0.8648 & 0.8712 & 0.8700 & 0.8636 & 0.8986 & 0.8688\\
& \multicolumn{1}{r}{$3$} & 0.0784 & 0.1984 & 0.3334 & 0.4464 & 0.5212 &
0.5798 & 0.6296 & 0.6602 & 0.6756 & 0.6708 & 0.6598 & 0.7104 & 0.6770\\
& \multicolumn{1}{r}{$2$} & 0.0176 & 0.0394 & 0.0786 & 0.1242 & 0.1672 &
0.2040 & 0.2462 & 0.2722 & 0.2912 & 0.2882 & 0.2810 & 0.3130 & 0.2764\\
& \multicolumn{1}{r}{$1/2$} & 0.3478 & 0.3584 & 0.3626 & 0.3482 & 0.3208 &
0.2840 & 0.2534 & 0.2102 & 0.1620 & 0.1020 & 0.0590 & 0.2280 & 0.2186\\
& \multicolumn{1}{r}{$1/3$} & 0.7302 & 0.7442 & 0.7490 & 0.7378 & 0.7086 &
0.6758 & 0.6302 & 0.5630 & 0.4662 & 0.3254 & 0.1820 & 0.5794 & 0.5436\\
& \multicolumn{1}{r}{$1/4$} & 0.9176 & 0.9216 & 0.9240 & 0.9170 & 0.9046 &
0.8882 & 0.8644 & 0.8230 & 0.7536 & 0.6104 & 0.4106 & 0.8304 & 0.7944\\
& \multicolumn{1}{r}{$1/5$} & 0.9768 & 0.9800 & 0.9804 & 0.9792 & 0.9744 &
0.9664 & 0.9568 & 0.9370 & 0.8954 & 0.8050 & 0.6226 & 0.9412 & 0.9184\\\hline
$100$ & \multicolumn{1}{r}{$5$} & 0.9924 & 0.9970 & 0.9984 & 0.9994 & 0.9996 &
1.0000 & 1.0000 & 1.0000 & 1.0000 & 1.0000 & 1.0000 & 1.0000 & 1.0000\\
& \multicolumn{1}{r}{$4$} & 0.9498 & 0.9694 & 0.9792 & 0.9882 & 0.9912 &
0.9930 & 0.9944 & 0.9950 & 0.9956 & 0.9956 & 0.9956 & 0.9960 & 0.9944\\
& \multicolumn{1}{r}{$3$} & 0.7272 & 0.8058 & 0.8526 & 0.8890 & 0.9118 &
0.9244 & 0.9344 & 0.9402 & 0.9448 & 0.9460 & 0.9496 & 0.9514 & 0.9432\\
& \multicolumn{1}{r}{$2$} & 0.1744 & 0.2498 & 0.3262 & 0.3994 & 0.4524 &
0.4938 & 0.5196 & 0.5410 & 0.5584 & 0.5664 & 0.5784 & 0.5700 & 0.5528\\
& \multicolumn{1}{r}{$1/2$} & 0.6496 & 0.6522 & 0.6446 & 0.6392 & 0.6246 &
0.6002 & 0.5552 & 0.5076 & 0.4480 & 0.3714 & 0.2976 & 0.4980 & 0.4886\\
& \multicolumn{1}{r}{$1/3$} & 0.9738 & 0.9744 & 0.9732 & 0.9724 & 0.9686 &
0.9632 & 0.9522 & 0.9374 & 0.9168 & 0.8776 & 0.8190 & 0.9322 & 0.9270\\
& \multicolumn{1}{r}{$1/4$} & 0.9996 & 0.9994 & 0.9996 & 0.9996 & 0.9994 &
0.9990 & 0.9982 & 0.9962 & 0.9948 & 0.9916 & 0.9842 & 0.9960 & 0.9952\\
& \multicolumn{1}{r}{$1/5$} & 1.0000 & 1.0000 & 1.0000 & 1.0000 & 1.0000 &
1.0000 & 1.0000 & 1.0000 & 1.0000 & 0.9994 & 0.9990 & 1.0000 & 0.9998\\\hline
$200$ & \multicolumn{1}{r}{$5$} & 1 & 1 & 1 & 1 & 1 & 1 & 1 & 1 & 1 & 1 & 1 &
1 & 1\\
& \multicolumn{1}{r}{$4$} & 1 & 1 & 1 & 1 & 1 & 1 & 1 & 1 & 1 & 1 & 1 & 1 &
1\\
& \multicolumn{1}{r}{$3$} & 0.9944 & 0.9960 & 0.9978 & 0.9978 & 0.9980 &
0.9986 & 0.9986 & 0.9992 & 0.9994 & 0.9994 & 0.9994 & 0.9994 & 0.9984\\
& \multicolumn{1}{r}{$2$} & 0.7046 & 0.7508 & 0.7844 & 0.8150 & 0.8328 &
0.8532 & 0.8636 & 0.8756 & 0.8828 & 0.8880 & 0.8892 & 0.8782 & 0.8738\\
& \multicolumn{1}{r}{$1/2$} & 0.9334 & 0.9314 & 0.9258 & 0.9190 & 0.9068 &
0.8976 & 0.8820 & 0.8634 & 0.8418 & 0.8122 & 0.7682 & 0.8432 & 0.8428\\
& \multicolumn{1}{r}{$1/3$} & 1.0000 & 1.0000 & 1.0000 & 1.0000 & 1.0000 &
1.0000 & 1.0000 & 0.9998 & 0.9996 & 0.9990 & 0.9984 & 0.9994 & 0.9994\\
& \multicolumn{1}{r}{$1/4$} & 1 & 1 & 1 & 1 & 1 & 1 & 1 & 1 & 1 & 1 & 1 & 1 &
1\\
& \multicolumn{1}{r}{$1/5$} & 1 & 1 & 1 & 1 & 1 & 1 & 1 & 1 & 1 & 1 & 1 & 1 &
1\\\hline
\end{tabular}
}
\end{center}
\end{table}

\begin{table}[ptb]
\caption{Empirical powers based on 5000 replications with simulated critical
values, when $\alpha=0.05$ and $\tau=0.3$.}
\begin{center}
{\tiny
\begin{tabular}
[c]{llccccccccccccc}
&  & \multicolumn{11}{c}{$^{0.05}T_{\phi_{\lambda}}^{(K)}$} & $\widetilde
{LRT}^{(K)}$ & $S^{(K)}$\\
&  & \multicolumn{11}{c}{$\lambda$} &  & \\
$K$ & \multicolumn{1}{r}{$\theta_{1}$} & $-1$ & $-0.9$ & $-0.8$ & $-0.7$ &
$-0.6$ & $-0.5$ & $-0.4$ & $-0.3$ & $-0.2$ & $-0.1$ & $0$ &  & \\\hline
$40$ & \multicolumn{1}{r}{$5$} & 0.2422 & 0.6320 & 0.8074 & 0.8870 & 0.9266 &
0.9446 & 0.9558 & 0.9594 & 0.9590 & 0.9552 & 0.9292 & 0.9714 & 0.9802\\
& \multicolumn{1}{r}{$4$} & 0.0998 & 0.3820 & 0.5938 & 0.7276 & 0.7992 &
0.8434 & 0.8698 & 0.8802 & 0.8798 & 0.8656 & 0.8056 & 0.9110 & 0.9316\\
& \multicolumn{1}{r}{$3$} & 0.0304 & 0.1362 & 0.2950 & 0.4262 & 0.5212 &
0.5876 & 0.6312 & 0.6526 & 0.6494 & 0.6310 & 0.5428 & 0.7146 & 0.7678\\
& \multicolumn{1}{r}{$2$} & 0.0154 & 0.0268 & 0.0596 & 0.1130 & 0.1640 &
0.2096 & 0.2480 & 0.2710 & 0.2698 & 0.2612 & 0.2050 & 0.3012 & 0.3548\\
& \multicolumn{1}{r}{$1/2$} & 0.2706 & 0.3078 & 0.3226 & 0.3210 & 0.3050 &
0.2738 & 0.2344 & 0.1768 & 0.1148 & 0.0638 & 0.0318 & 0.2558 & 0.3336\\
& \multicolumn{1}{r}{$1/3$} & 0.6232 & 0.6788 & 0.7024 & 0.7090 & 0.6916 &
0.6582 & 0.6106 & 0.5290 & 0.3984 & 0.2376 & 0.0612 & 0.6322 & 0.7250\\
& \multicolumn{1}{r}{$1/4$} & 0.8490 & 0.8866 & 0.9010 & 0.9034 & 0.8952 &
0.8744 & 0.8436 & 0.7818 & 0.6736 & 0.4666 & 0.1452 & 0.8522 & 0.9070\\
& \multicolumn{1}{r}{$1/5$} & 0.9554 & 0.9686 & 0.9726 & 0.9748 & 0.9706 &
0.9630 & 0.9538 & 0.9304 & 0.8640 & 0.7168 & 0.3112 & 0.9578 & 0.9790\\\hline
$50$ & \multicolumn{1}{r}{$5$} & 0.7892 & 0.9106 & 0.9534 & 0.9732 & 0.9846 &
0.9876 & 0.9900 & 0.9916 & 0.9924 & 0.9914 & 0.9906 & 0.9926 & 0.9942\\
& \multicolumn{1}{r}{$4$} & 0.5456 & 0.7420 & 0.8446 & 0.8902 & 0.9214 &
0.9394 & 0.9492 & 0.9562 & 0.9586 & 0.9548 & 0.9476 & 0.9600 & 0.9716\\
& \multicolumn{1}{r}{$3$} & 0.2176 & 0.4054 & 0.5546 & 0.6530 & 0.7130 &
0.7602 & 0.7944 & 0.8128 & 0.8210 & 0.8092 & 0.7906 & 0.8260 & 0.8620\\
& \multicolumn{1}{r}{$2$} & 0.0290 & 0.0754 & 0.1434 & 0.2038 & 0.2596 &
0.3056 & 0.3486 & 0.3726 & 0.3862 & 0.3756 & 0.3548 & 0.3872 & 0.4368\\
& \multicolumn{1}{r}{$1/2$} & 0.4212 & 0.4408 & 0.4540 & 0.4442 & 0.4168 &
0.3792 & 0.3452 & 0.2934 & 0.2240 & 0.1364 & 0.0692 & 0.3404 & 0.4220\\
& \multicolumn{1}{r}{$1/3$} & 0.8252 & 0.8420 & 0.8488 & 0.8474 & 0.8338 &
0.8108 & 0.7832 & 0.7352 & 0.6480 & 0.5084 & 0.3080 & 0.7702 & 0.8388\\
& \multicolumn{1}{r}{$1/4$} & 0.9650 & 0.9694 & 0.9722 & 0.9710 & 0.9676 &
0.9626 & 0.9546 & 0.9338 & 0.8990 & 0.8160 & 0.6428 & 0.9474 & 0.9710\\
& \multicolumn{1}{r}{$1/5$} & 0.9946 & 0.9956 & 0.9956 & 0.9954 & 0.9952 &
0.9940 & 0.9918 & 0.9880 & 0.9778 & 0.9436 & 0.8558 & 0.9906 & 0.9952\\\hline
$100$ & \multicolumn{1}{r}{$5$} & 0.9998 & 1.0000 & 1.0000 & 1.0000 & 1.0000 &
1.0000 & 1.0000 & 1.0000 & 1.0000 & 1.0000 & 1.0000 & 1.0000 & 1.0000\\
& \multicolumn{1}{r}{$4$} & 0.9944 & 0.9978 & 0.9984 & 0.9990 & 0.9994 &
0.9994 & 0.9996 & 0.9996 & 0.9996 & 0.9996 & 0.9996 & 0.9996 & 1.0000\\
& \multicolumn{1}{r}{$3$} & 0.9096 & 0.9414 & 0.9608 & 0.9746 & 0.9808 &
0.9842 & 0.9862 & 0.9874 & 0.9890 & 0.9892 & 0.9888 & 0.9882 & 0.9938\\
& \multicolumn{1}{r}{$2$} & 0.3412 & 0.4398 & 0.5190 & 0.5860 & 0.6388 &
0.6734 & 0.6914 & 0.7066 & 0.7196 & 0.7208 & 0.7222 & 0.7082 & 0.7864\\
& \multicolumn{1}{r}{$1/2$} & 0.7582 & 0.7638 & 0.7598 & 0.7564 & 0.7458 &
0.7266 & 0.6898 & 0.6508 & 0.5962 & 0.5188 & 0.4360 & 0.6646 & 0.7658\\
& \multicolumn{1}{r}{$1/3$} & 0.9928 & 0.9932 & 0.9932 & 0.9932 & 0.9928 &
0.9924 & 0.9894 & 0.9850 & 0.9770 & 0.9638 & 0.9398 & 0.9876 & 0.9956\\
& \multicolumn{1}{r}{$1/4$} & 1.0000 & 1.0000 & 1.0000 & 1.0000 & 1.0000 &
1.0000 & 0.9998 & 0.9998 & 0.9996 & 0.9990 & 0.9984 & 0.9998 & 1.0000\\
& \multicolumn{1}{r}{$1/5$} & 1.0000 & 1.0000 & 1.0000 & 1.0000 & 1.0000 &
1.0000 & 1.0000 & 1.0000 & 1.0000 & 1.0000 & 0.9998 & 1.0000 & 1.0000\\\hline
$200$ & \multicolumn{1}{r}{$5$} & 1.0000 & 1.0000 & 1.0000 & 1.0000 & 1.0000 &
1.0000 & 1.0000 & 1.0000 & 1.0000 & 1.0000 & 1.0000 & 1.0000 & 1.0000\\
& \multicolumn{1}{r}{$4$} & 1.0000 & 1.0000 & 1.0000 & 1.0000 & 1.0000 &
1.0000 & 1.0000 & 1.0000 & 1.0000 & 1.0000 & 1.0000 & 1.0000 & 1.0000\\
& \multicolumn{1}{r}{$3$} & 1.0000 & 1.0000 & 1.0000 & 1.0000 & 1.0000 &
1.0000 & 1.0000 & 1.0000 & 1.0000 & 1.0000 & 1.0000 & 1.0000 & 1.0000\\
& \multicolumn{1}{r}{$2$} & 0.8920 & 0.9120 & 0.9256 & 0.9394 & 0.9464 &
0.9532 & 0.9558 & 0.9596 & 0.9618 & 0.9634 & 0.9634 & 0.9542 & 0.9748\\
& \multicolumn{1}{r}{$1/2$} & 0.9762 & 0.9768 & 0.9744 & 0.9716 & 0.9658 &
0.9628 & 0.9576 & 0.9492 & 0.9372 & 0.9214 & 0.8970 & 0.9442 & 0.9796\\
& \multicolumn{1}{r}{$1/3$} & 1.0000 & 1.0000 & 1.0000 & 1.0000 & 1.0000 &
1.0000 & 1.0000 & 1.0000 & 1.0000 & 1.0000 & 1.0000 & 1.0000 & 1.0000\\
& \multicolumn{1}{r}{$1/4$} & 1.0000 & 1.0000 & 1.0000 & 1.0000 & 1.0000 &
1.0000 & 1.0000 & 1.0000 & 1.0000 & 1.0000 & 1.0000 & 1.0000 & 1.0000\\
& \multicolumn{1}{r}{$1/5$} & 1.0000 & 1.0000 & 1.0000 & 1.0000 & 1.0000 &
1.0000 & 1.0000 & 1.0000 & 1.0000 & 1.0000 & 1.0000 & 1.0000 & 1.0000\\\hline
\end{tabular}
}
\end{center}
\end{table}

\begin{table}[ptb]
\caption{Empirical powers based on 5000 replications with simulated critical
values, when $\alpha=0.05$ and $\tau=0.5$.}
\begin{center}
{\tiny
\begin{tabular}
[c]{llccccccccccccc}
&  & \multicolumn{11}{c}{$^{0.05}T_{\phi_{\lambda}}^{(K)}$} & $\widetilde
{LRT}^{(K)}$ & $S^{(K)}$\\
&  & \multicolumn{11}{c}{$\lambda$} &  & \\
$K$ & \multicolumn{1}{r}{$\theta_{1}$} & $-1$ & $-0.9$ & $-0.8$ & $-0.7$ &
$-0.6$ & $-0.5$ & $-0.4$ & $-0.3$ & $-0.2$ & $-0.1$ & $0$ &  & \\\hline
$40$ & \multicolumn{1}{r}{$5$} & 0.4646 & 0.8180 & 0.9282 & 0.9646 & 0.9794 &
0.9856 & 0.9886 & 0.9898 & 0.9888 & 0.9870 & 0.9752 & 0.9888 & 0.9970\\
& \multicolumn{1}{r}{$4$} & 0.2220 & 0.5834 & 0.7678 & 0.8668 & 0.9106 &
0.9346 & 0.9486 & 0.9524 & 0.9508 & 0.9420 & 0.9016 & 0.9492 & 0.9784\\
& \multicolumn{1}{r}{$3$} & 0.0632 & 0.2624 & 0.4530 & 0.5914 & 0.6792 &
0.7304 & 0.7654 & 0.7794 & 0.7726 & 0.7496 & 0.6620 & 0.7678 & 0.8726\\
& \multicolumn{1}{r}{$2$} & 0.0262 & 0.0512 & 0.1070 & 0.1850 & 0.2406 &
0.2910 & 0.3372 & 0.3538 & 0.3492 & 0.3302 & 0.2638 & 0.3314 & 0.4716\\
& \multicolumn{1}{r}{$1/2$} & 0.2752 & 0.3292 & 0.3488 & 0.3566 & 0.3402 &
0.3112 & 0.2728 & 0.2056 & 0.1296 & 0.0630 & 0.0234 & 0.3522 & 0.4816\\
& \multicolumn{1}{r}{$1/3$} & 0.6534 & 0.7274 & 0.7584 & 0.7644 & 0.7544 &
0.7248 & 0.6762 & 0.5874 & 0.4666 & 0.2816 & 0.0596 & 0.7614 & 0.8612\\
& \multicolumn{1}{r}{$1/4$} & 0.8972 & 0.9270 & 0.9368 & 0.9408 & 0.9366 &
0.9264 & 0.9098 & 0.8688 & 0.7806 & 0.6120 & 0.2294 & 0.9372 & 0.9720\\
& \multicolumn{1}{r}{$1/5$} & 0.9712 & 0.9834 & 0.9872 & 0.9882 & 0.9872 &
0.9846 & 0.9790 & 0.9632 & 0.9298 & 0.8374 & 0.4684 & 0.9890 & 0.9970\\\hline
$50$ & \multicolumn{1}{r}{$5$} & 0.9212 & 0.9748 & 0.9882 & 0.9950 & 0.9970 &
0.9976 & 0.9980 & 0.9980 & 0.9982 & 0.9980 & 0.9974 & 0.9980 & 0.9998\\
& \multicolumn{1}{r}{$4$} & 0.7448 & 0.8846 & 0.9438 & 0.9652 & 0.9760 &
0.9810 & 0.9840 & 0.9858 & 0.9866 & 0.9854 & 0.9830 & 0.9808 & 0.9928\\
& \multicolumn{1}{r}{$3$} & 0.3800 & 0.5908 & 0.7284 & 0.8006 & 0.8438 &
0.8680 & 0.8902 & 0.9006 & 0.9036 & 0.8950 & 0.8770 & 0.8716 & 0.9408\\
& \multicolumn{1}{r}{$2$} & 0.0620 & 0.1470 & 0.2340 & 0.3122 & 0.3676 &
0.4146 & 0.4596 & 0.4844 & 0.4928 & 0.4766 & 0.4504 & 0.4298 & 0.5818\\
& \multicolumn{1}{r}{$1/2$} & 0.4470 & 0.4760 & 0.4904 & 0.4836 & 0.4628 &
0.4352 & 0.3970 & 0.3378 & 0.2646 & 0.1606 & 0.0706 & 0.4448 & 0.5832\\
& \multicolumn{1}{r}{$1/3$} & 0.8734 & 0.8894 & 0.9010 & 0.8998 & 0.8882 &
0.8728 & 0.8534 & 0.8130 & 0.7436 & 0.6166 & 0.4168 & 0.8814 & 0.9432\\
& \multicolumn{1}{r}{$1/4$} & 0.9798 & 0.9834 & 0.9850 & 0.9850 & 0.9834 &
0.9814 & 0.9762 & 0.9652 & 0.9444 & 0.8972 & 0.7766 & 0.9816 & 0.9948\\
& \multicolumn{1}{r}{$1/5$} & 0.9978 & 0.9986 & 0.9986 & 0.9986 & 0.9986 &
0.9976 & 0.9970 & 0.9948 & 0.9888 & 0.9794 & 0.9396 & 0.9974 & 0.9994\\\hline
$100$ & \multicolumn{1}{r}{$5$} & 1 & 1 & 1 & 1 & 1 & 1 & 1 & 1 & 1 & 1 & 1 &
1 & 1\\
& \multicolumn{1}{r}{$4$} & 0.9996 & 0.9996 & 1.0000 & 1.0000 & 1.0000 &
1.0000 & 1.0000 & 1.0000 & 1.0000 & 1.0000 & 1.0000 & 1.0000 & 1.0000\\
& \multicolumn{1}{r}{$3$} & 0.9730 & 0.9840 & 0.9900 & 0.9930 & 0.9950 &
0.9964 & 0.9966 & 0.9970 & 0.9972 & 0.9970 & 0.9966 & 0.9944 & 0.9988\\
& \multicolumn{1}{r}{$2$} & 0.4946 & 0.5952 & 0.6668 & 0.7292 & 0.7714 &
0.7966 & 0.8112 & 0.8190 & 0.8260 & 0.8246 & 0.8246 & 0.7812 & 0.8872\\
& \multicolumn{1}{r}{$1/2$} & 0.8016 & 0.8096 & 0.8104 & 0.8090 & 0.7988 &
0.7826 & 0.7532 & 0.7186 & 0.6708 & 0.5982 & 0.5200 & 0.7714 & 0.8866\\
& \multicolumn{1}{r}{$1/3$} & 0.9978 & 0.9980 & 0.9980 & 0.9980 & 0.9976 &
0.9974 & 0.9966 & 0.9954 & 0.9930 & 0.9880 & 0.9784 & 0.9970 & 0.9994\\
& \multicolumn{1}{r}{$1/4$} & 1.0000 & 1.0000 & 1.0000 & 1.0000 & 1.0000 &
1.0000 & 1.0000 & 1.0000 & 1.0000 & 0.9998 & 0.9996 & 1.0000 & 1.0000\\
& \multicolumn{1}{r}{$1/5$} & 1 & 1 & 1 & 1 & 1 & 1 & 1 & 1 & 1 & 1 & 1 & 1 &
1\\\hline
$200$ & \multicolumn{1}{r}{$5$} & 1 & 1 & 1 & 1 & 1 & 1 & 1 & 1 & 1 & 1 & 1 &
1 & 1\\
& \multicolumn{1}{r}{$4$} & 1 & 1 & 1 & 1 & 1 & 1 & 1 & 1 & 1 & 1 & 1 & 1 &
1\\
& \multicolumn{1}{r}{$3$} & 1 & 1 & 1 & 1 & 1 & 1 & 1 & 1 & 1 & 1 & 1 & 1 &
1\\
& \multicolumn{1}{r}{$2$} & 0.9562 & 0.9678 & 0.9736 & 0.9776 & 0.9808 &
0.9848 & 0.9860 & 0.9878 & 0.9882 & 0.9890 & 0.9884 & 0.9788 & 0.9956\\
& \multicolumn{1}{r}{$1/2$} & 0.9858 & 0.9858 & 0.9856 & 0.9852 & 0.9832 &
0.9824 & 0.9784 & 0.9740 & 0.9668 & 0.9596 & 0.9466 & 0.9776 & 0.9952\\
& \multicolumn{1}{r}{$1/3$} & 1 & 1 & 1 & 1 & 1 & 1 & 1 & 1 & 1 & 1 & 1 & 1 &
1\\
& \multicolumn{1}{r}{$1/4$} & 1 & 1 & 1 & 1 & 1 & 1 & 1 & 1 & 1 & 1 & 1 & 1 &
1\\
& \multicolumn{1}{r}{$1/5$} & 1 & 1 & 1 & 1 & 1 & 1 & 1 & 1 & 1 & 1 & 1 & 1 &
1\\\hline
\end{tabular}
}
\end{center}
\end{table}

\begin{table}[ptb]
\caption{Empirical powers based on 5000 replications with simulated critical
values, when $\alpha=0.01$ and $\tau=0.2$.}
\begin{center}
{\tiny
\begin{tabular}
[c]{llccccccccccccc}
&  & \multicolumn{11}{c}{$^{0.05}T_{\phi_{\lambda}}^{(K)}$} & $\widetilde
{LRT}^{(K)}$ & $S^{(K)}$\\
&  & \multicolumn{11}{c}{$\lambda$} &  & \\
$K$ & \multicolumn{1}{r}{$\theta_{1}$} & $-1$ & $-0.9$ & $-0.8$ & $-0.7$ &
$-0.6$ & $-0.5$ & $-0.4$ & $-0.3$ & $-0.2$ & $-0.1$ & $0$ &  & \\\hline
$40$ & \multicolumn{1}{r}{$5$} & 0.0008 & 0.0052 & 0.1090 & 0.3356 & 0.5070 &
0.6232 & 0.6720 & 0.6876 & 0.6522 & 0.5440 & 0.3070 & 0.8272 & 0.7804\\
& \multicolumn{1}{r}{$4$} & 0.0012 & 0.0012 & 0.0376 & 0.1642 & 0.3068 &
0.4186 & 0.4670 & 0.4846 & 0.4460 & 0.3378 & 0.1530 & 0.6682 & 0.6016\\
& \multicolumn{1}{r}{$3$} & 0.0024 & 0.0016 & 0.0096 & 0.0498 & 0.1130 &
0.1808 & 0.2172 & 0.2296 & 0.2008 & 0.1314 & 0.0526 & 0.3866 & 0.3292\\
& \multicolumn{1}{r}{$2$} & 0.0016 & 0.0012 & 0.0012 & 0.0052 & 0.0172 &
0.0332 & 0.0440 & 0.0474 & 0.0410 & 0.0252 & 0.0166 & 0.1028 & 0.0790\\
& \multicolumn{1}{r}{$1/2$} & 0.0374 & 0.0540 & 0.0752 & 0.0864 & 0.0854 &
0.0714 & 0.0500 & 0.0292 & 0.0122 & 0.0076 & 0.0084 & 0.0520 & 0.0562\\
& \multicolumn{1}{r}{$1/3$} & 0.1066 & 0.1980 & 0.2646 & 0.3010 & 0.2970 &
0.2672 & 0.1954 & 0.1084 & 0.0340 & 0.0062 & 0.0064 & 0.1964 & 0.1628\\
& \multicolumn{1}{r}{$1/4$} & 0.2234 & 0.3960 & 0.4962 & 0.5392 & 0.5340 &
0.4936 & 0.3946 & 0.2544 & 0.0916 & 0.0078 & 0.0062 & 0.3790 & 0.3128\\
& \multicolumn{1}{r}{$1/5$} & 0.3792 & 0.6034 & 0.7102 & 0.7444 & 0.7404 &
0.7098 & 0.6128 & 0.4410 & 0.2002 & 0.0200 & 0.0064 & 0.5898 & 0.4872\\\hline
$50$ & \multicolumn{1}{r}{$5$} & 0.0074 & 0.1362 & 0.4184 & 0.6102 & 0.7354 &
0.7920 & 0.8300 & 0.8444 & 0.8328 & 0.8062 & 0.7132 & 0.9172 & 0.8758\\
& \multicolumn{1}{r}{$4$} & 0.0012 & 0.0378 & 0.1942 & 0.3598 & 0.4992 &
0.5872 & 0.6448 & 0.6678 & 0.6464 & 0.6058 & 0.4746 & 0.7866 & 0.7118\\
& \multicolumn{1}{r}{$3$} & 0.0004 & 0.0066 & 0.0484 & 0.1280 & 0.2168 &
0.2868 & 0.3448 & 0.3728 & 0.3454 & 0.3044 & 0.2002 & 0.5214 & 0.4384\\
& \multicolumn{1}{r}{$2$} & 0.0022 & 0.0018 & 0.0060 & 0.0172 & 0.0360 &
0.0566 & 0.0776 & 0.0848 & 0.0788 & 0.0652 & 0.0342 & 0.1514 & 0.1132\\
& \multicolumn{1}{r}{$1/2$} & 0.1082 & 0.1234 & 0.1442 & 0.1416 & 0.1332 &
0.1082 & 0.0796 & 0.0468 & 0.0170 & 0.0086 & 0.0064 & 0.0746 & 0.0678\\
& \multicolumn{1}{r}{$1/3$} & 0.3668 & 0.4150 & 0.4554 & 0.4554 & 0.4396 &
0.3884 & 0.3212 & 0.2232 & 0.1006 & 0.0236 & 0.0046 & 0.2974 & 0.2274\\
& \multicolumn{1}{r}{$1/4$} & 0.6416 & 0.7016 & 0.7418 & 0.7456 & 0.7302 &
0.6830 & 0.6080 & 0.4960 & 0.2808 & 0.0830 & 0.0078 & 0.5726 & 0.4620\\
& \multicolumn{1}{r}{$1/5$} & 0.8276 & 0.8642 & 0.8902 & 0.8934 & 0.8834 &
0.8550 & 0.8086 & 0.7010 & 0.4772 & 0.1936 & 0.0072 & 0.7738 & 0.6572\\\hline
$100$ & \multicolumn{1}{r}{$5$} & 0.8036 & 0.9376 & 0.9370 & 0.9886 & 0.9944 &
0.9968 & 0.9974 & 0.9976 & 0.9980 & 0.9972 & 0.9970 & 0.9994 & 0.9982\\
& \multicolumn{1}{r}{$4$} & 0.4946 & 0.7542 & 0.7508 & 0.9292 & 0.9576 &
0.9670 & 0.9736 & 0.9762 & 0.9774 & 0.9746 & 0.9682 & 0.9880 & 0.9810\\
& \multicolumn{1}{r}{$3$} & 0.1256 & 0.3418 & 0.3380 & 0.6714 & 0.7520 &
0.7930 & 0.8222 & 0.8354 & 0.8434 & 0.8282 & 0.8006 & 0.8926 & 0.8592\\
& \multicolumn{1}{r}{$2$} & 0.0028 & 0.0206 & 0.0200 & 0.1312 & 0.1924 &
0.2282 & 0.2700 & 0.2958 & 0.3078 & 0.2834 & 0.2438 & 0.3992 & 0.3450\\
& \multicolumn{1}{r}{$1/2$} & 0.3122 & 0.3494 & 0.2528 & 0.3782 & 0.3740 &
0.3360 & 0.2902 & 0.2300 & 0.1668 & 0.0898 & 0.0276 & 0.2876 & 0.2634\\
& \multicolumn{1}{r}{$1/3$} & 0.8326 & 0.8602 & 0.7892 & 0.8864 & 0.8838 &
0.8550 & 0.8232 & 0.7646 & 0.6772 & 0.5060 & 0.2610 & 0.8132 & 0.7532\\
& \multicolumn{1}{r}{$1/4$} & 0.9864 & 0.9898 & 0.9786 & 0.9930 & 0.9928 &
0.9906 & 0.9858 & 0.9728 & 0.9526 & 0.8746 & 0.6808 & 0.9830 & 0.9654\\
& \multicolumn{1}{r}{$1/5$} & 0.9990 & 0.9992 & 0.9988 & 0.9994 & 0.9994 &
0.9992 & 0.9990 & 0.9986 & 0.9956 & 0.9816 & 0.9070 & 0.9990 & 0.9972\\\hline
$200$ & \multicolumn{1}{r}{$5$} & 1 & 1 & 1 & 1 & 1 & 1 & 1 & 1 & 1 & 1 & 1 &
1 & 1\\
& \multicolumn{1}{r}{$4$} & 0.9986 & 0.9996 & 0.9998 & 0.9998 & 0.9998 &
1.0000 & 1.0000 & 1.0000 & 1.0000 & 1.0000 & 1.0000 & 1.0000 & 0.9998\\
& \multicolumn{1}{r}{$3$} & 0.9506 & 0.9726 & 0.9842 & 0.9890 & 0.9922 &
0.9930 & 0.9936 & 0.9936 & 0.9940 & 0.9940 & 0.9936 & 0.9958 & 0.9940\\
& \multicolumn{1}{r}{$2$} & 0.3052 & 0.4120 & 0.4978 & 0.5508 & 0.5938 &
0.6486 & 0.6676 & 0.6688 & 0.6822 & 0.6784 & 0.6614 & 0.7446 & 0.6786\\
& \multicolumn{1}{r}{$1/2$} & 0.7840 & 0.7912 & 0.7906 & 0.7706 & 0.7428 &
0.7314 & 0.6872 & 0.6104 & 0.5478 & 0.4524 & 0.3242 & 0.6528 & 0.5752\\
& \multicolumn{1}{r}{$1/3$} & 0.9986 & 0.9988 & 0.9986 & 0.9984 & 0.9980 &
0.9978 & 0.9964 & 0.9930 & 0.9906 & 0.9818 & 0.9588 & 0.9948 & 0.9908\\
& \multicolumn{1}{r}{$1/4$} & 1 & 1 & 1 & 1 & 1 & 1 & 1 & 1 & 1 & 1 & 1 & 1 &
1\\
& \multicolumn{1}{r}{$1/5$} & 1 & 1 & 1 & 1 & 1 & 1 & 1 & 1 & 1 & 1 & 1 & 1 &
1\\\hline
\end{tabular}
}
\end{center}
\end{table}\begin{table}[ptb]
\caption{Empirical powers based on 5000 replications with simulated critical
values, when $\alpha=0.01$ and $\tau=0.3$.}
\begin{center}
{\tiny
\begin{tabular}
[c]{llccccccccccccc}
&  & \multicolumn{11}{c}{$^{0.05}T_{\phi_{\lambda}}^{(K)}$} & $\widetilde
{LRT}^{(K)}$ & $S^{(K)}$\\
&  & \multicolumn{11}{c}{$\lambda$} &  & \\
$K$ & \multicolumn{1}{r}{$\theta_{1}$} & $-1$ & $-0.9$ & $-0.8$ & $-0.7$ &
$-0.6$ & $-0.5$ & $-0.4$ & $-0.3$ & $-0.2$ & $-0.1$ & $0$ &  & \\\hline
$40$ & \multicolumn{1}{r}{$5$} & 0.0024 & 0.0250 & 0.3266 & 0.6294 & 0.7738 &
0.8456 & 0.8696 & 0.8728 & 0.8476 & 0.7584 & 0.4934 & 0.9198 & 0.9390\\
& \multicolumn{1}{r}{$4$} & 0.0010 & 0.0066 & 0.1342 & 0.3792 & 0.5456 &
0.6534 & 0.6912 & 0.6982 & 0.6552 & 0.5286 & 0.2596 & 0.7906 & 0.8322\\
& \multicolumn{1}{r}{$3$} & 0.0014 & 0.0016 & 0.0330 & 0.1284 & 0.2496 &
0.3452 & 0.3894 & 0.3972 & 0.3520 & 0.2400 & 0.0810 & 0.5052 & 0.5582\\
& \multicolumn{1}{r}{$2$} & 0.0030 & 0.0018 & 0.0032 & 0.0164 & 0.0390 &
0.0706 & 0.0874 & 0.0930 & 0.0778 & 0.0426 & 0.0220 & 0.1404 & 0.1742\\
& \multicolumn{1}{r}{$1/2$} & 0.0362 & 0.0698 & 0.1050 & 0.1244 & 0.1268 &
0.1134 & 0.0782 & 0.0390 & 0.0116 & 0.0044 & 0.0060 & 0.0928 & 0.1324\\
& \multicolumn{1}{r}{$1/3$} & 0.1222 & 0.2812 & 0.3958 & 0.4502 & 0.4544 &
0.4240 & 0.3390 & 0.2250 & 0.0788 & 0.0068 & 0.0066 & 0.3622 & 0.4640\\
& \multicolumn{1}{r}{$1/4$} & 0.2712 & 0.5376 & 0.6680 & 0.7186 & 0.7256 &
0.7028 & 0.6168 & 0.4604 & 0.2142 & 0.0240 & 0.0076 & 0.6408 & 0.7270\\
& \multicolumn{1}{r}{$1/5$} & 0.4944 & 0.7690 & 0.8636 & 0.8964 & 0.8992 &
0.8834 & 0.8342 & 0.7182 & 0.4458 & 0.0680 & 0.0062 & 0.8458 & 0.8972\\\hline
$50$ & \multicolumn{1}{r}{$5$} & 0.0464 & 0.4022 & 0.7184 & 0.8592 & 0.9172 &
0.9410 & 0.9532 & 0.9564 & 0.9504 & 0.9346 & 0.8728 & 0.9726 & 0.9794\\
& \multicolumn{1}{r}{$4$} & 0.0096 & 0.1638 & 0.4522 & 0.6352 & 0.7558 &
0.8114 & 0.8442 & 0.8552 & 0.8360 & 0.7942 & 0.6636 & 0.8898 & 0.9100\\
& \multicolumn{1}{r}{$3$} & 0.0024 & 0.0306 & 0.1540 & 0.2950 & 0.4236 &
0.5010 & 0.5522 & 0.5706 & 0.5384 & 0.4828 & 0.3334 & 0.6454 & 0.6816\\
& \multicolumn{1}{r}{$2$} & 0.0022 & 0.0018 & 0.0060 & 0.0172 & 0.0360 &
0.0566 & 0.0776 & 0.0848 & 0.0788 & 0.0652 & 0.0342 & 0.1514 & 0.1132\\
& \multicolumn{1}{r}{$1/2$} & 0.1216 & 0.1574 & 0.1938 & 0.1978 & 0.1908 &
0.1634 & 0.1244 & 0.0778 & 0.0284 & 0.0080 & 0.0028 & 0.1340 & 0.1830\\
& \multicolumn{1}{r}{$1/3$} & 0.4600 & 0.5436 & 0.6072 & 0.6216 & 0.6144 &
0.5704 & 0.5034 & 0.3912 & 0.2112 & 0.0578 & 0.0054 & 0.5132 & 0.5944\\
& \multicolumn{1}{r}{$1/4$} & 0.7796 & 0.8450 & 0.8824 & 0.8894 & 0.8856 &
0.8640 & 0.8204 & 0.7320 & 0.5192 & 0.2248 & 0.0104 & 0.8258 & 0.8724\\
& \multicolumn{1}{r}{$1/5$} & 0.9242 & 0.9548 & 0.9686 & 0.9720 & 0.9714 &
0.9628 & 0.9434 & 0.9034 & 0.7804 & 0.4776 & 0.0356 & 0.9452 & 0.9638\\\hline
$100$ & \multicolumn{1}{r}{$5$} & 0.9700 & 0.9946 & 0.9942 & 0.9996 & 1.0000 &
1.0000 & 1.0000 & 1.0000 & 1.0000 & 1.0000 & 1.0000 & 1.0000 & 1\\
& \multicolumn{1}{r}{$4$} & 0.8182 & 0.9414 & 0.9410 & 0.9894 & 0.9952 &
0.9968 & 0.9978 & 0.9980 & 0.9980 & 0.9976 & 0.9956 & 0.9990 & 0.9992\\
& \multicolumn{1}{r}{$3$} & 0.3640 & 0.6342 & 0.6304 & 0.8772 & 0.9188 &
0.9316 & 0.9434 & 0.9470 & 0.9478 & 0.9388 & 0.9228 & 0.9662 & 0.9768\\
& \multicolumn{1}{r}{$2$} & 0.0168 & 0.0760 & 0.0742 & 0.2708 & 0.3620 &
0.4082 & 0.4462 & 0.4660 & 0.4732 & 0.4390 & 0.3826 & 0.5322 & 0.5996\\
& \multicolumn{1}{r}{$1/2$} & 0.4036 & 0.4582 & 0.3680 & 0.5060 & 0.5096 &
0.4748 & 0.4270 & 0.3678 & 0.2806 & 0.1554 & 0.0584 & 0.4560 & 0.5560\\
& \multicolumn{1}{r}{$1/3$} & 0.9244 & 0.9478 & 0.9132 & 0.9624 & 0.9634 &
0.9534 & 0.9404 & 0.9134 & 0.8720 & 0.7588 & 0.5352 & 0.9472 & 0.9708\\
& \multicolumn{1}{r}{$1/4$} & 0.9978 & 0.9984 & 0.9972 & 0.9990 & 0.9990 &
0.9986 & 0.9984 & 0.9972 & 0.9948 & 0.9774 & 0.9190 & 0.9986 & \\
& \multicolumn{1}{r}{$1/5$} & 0.9998 & 1.0000 & 0.9998 & 1.0000 & 1.0000 &
1.0000 & 0.9998 & 0.9998 & 0.9998 & 0.9990 & 0.9936 & 0.9998 & 1\\\hline
$200$ & \multicolumn{1}{r}{$5$} & 1 & 1 & 1 & 1 & 1 & 1 & 1 & 1 & 1 & 1 & 1 &
1 & 1\\
& \multicolumn{1}{r}{$4$} & 1 & 1 & 1 & 1 & 1 & 1 & 1 & 1 & 1 & 1 & 1 & 1 &
1\\
& \multicolumn{1}{r}{$3$} & 0.9962 & 0.9982 & 0.9992 & 0.9996 & 0.9998 &
0.9998 & 0.9998 & 0.9998 & 1.0000 & 1.0000 & 1.0000 & 0.9998 & 0.9998\\
& \multicolumn{1}{r}{$2$} & 0.5808 & 0.6866 & 0.7604 & 0.7912 & 0.8204 &
0.8556 & 0.8674 & 0.8650 & 0.8700 & 0.8632 & 0.8446 & 0.8886 & 0.9178\\
& \multicolumn{1}{r}{$1/2$} & 0.8892 & 0.8958 & 0.9004 & 0.8890 & 0.8772 &
0.8700 & 0.8424 & 0.7958 & 0.7442 & 0.6596 & 0.5364 & 0.8358 & 0.8964\\
& \multicolumn{1}{r}{$1/3$} & 1.0000 & 1.0000 & 1.0000 & 1.0000 & 1.0000 &
1.0000 & 1.0000 & 1.0000 & 0.9998 & 0.9986 & 0.9950 & 1.0000 & 1\\
& \multicolumn{1}{r}{$1/4$} & 1.0000 & 1.0000 & 1.0000 & 1.0000 & 1.0000 &
1.0000 & 1.0000 & 1.0000 & 0.9998 & 0.9986 & 0.9950 & 1.0000 & 1\\
& \multicolumn{1}{r}{$1/5$} & 1 & 1 & 1 & 1 & 1 & 1 & 1 & 1 & 1 & 1 & 1 & 1 &
1\\\hline
\end{tabular}
}
\end{center}
\end{table}

\begin{table}[ptb]
\caption{Empirical powers based on 5000 replications with simulated critical
values, when $\alpha=0.01$ and $\tau=0.5$.}
\begin{center}
{\tiny
\begin{tabular}
[c]{llccccccccccccc}
&  & \multicolumn{11}{c}{$^{0.05}T_{\phi_{\lambda}}^{(K)}$} & $\widetilde
{LRT}^{(K)}$ & $S^{(K)}$\\
&  & \multicolumn{11}{c}{$\lambda$} &  & \\
$K$ & \multicolumn{1}{r}{$\theta_{1}$} & $-1$ & $-0.9$ & $-0.8$ & $-0.7$ &
$-0.6$ & $-0.5$ & $-0.4$ & $-0.3$ & $-0.2$ & $-0.1$ & $0$ &  & \\\hline
$40$ & \multicolumn{1}{r}{$5$} & 0.0038 & 0.0902 & 0.5536 & 0.8170 & 0.9100 &
0.9458 & 0.9562 & 0.9562 & 0.9412 & 0.8818 & 0.6374 & 0.9538 & 0.9834\\
& \multicolumn{1}{r}{$4$} & 0.0022 & 0.0242 & 0.2914 & 0.5764 & 0.7260 &
0.8120 & 0.8376 & 0.8388 & 0.7998 & 0.6846 & 0.3880 & 0.8344 & 0.9252\\
& \multicolumn{1}{r}{$3$} & 0.0028 & 0.0046 & 0.0790 & 0.2530 & 0.4036 &
0.5014 & 0.5406 & 0.5434 & 0.4888 & 0.3646 & 0.1450 & 0.5366 & 0.7054\\
& \multicolumn{1}{r}{$2$} & 0.0038 & 0.0028 & 0.0080 & 0.0352 & 0.0780 &
0.1204 & 0.1468 & 0.1502 & 0.1238 & 0.0758 & 0.0346 & 0.1434 & 0.2444\\
& \multicolumn{1}{r}{$1/2$} & 0.0294 & 0.0674 & 0.1164 & 0.1470 & 0.1508 &
0.1376 & 0.0956 & 0.0486 & 0.0124 & 0.0036 & 0.0056 & 0.1586 & 0.2562\\
& \multicolumn{1}{r}{$1/3$} & 0.1056 & 0.3064 & 0.4406 & 0.5010 & 0.5104 &
0.4900 & 0.4144 & 0.2788 & 0.1024 & 0.0074 & 0.0050 & 0.5220 & 0.6962\\
& \multicolumn{1}{r}{$1/4$} & 0.3086 & 0.6340 & 0.7632 & 0.8088 & 0.8198 &
0.7996 & 0.7414 & 0.6150 & 0.3556 & 0.0442 & 0.0044 & 0.8286 & 0.9176\\
& \multicolumn{1}{r}{$1/5$} & 0.5660 & 0.8496 & 0.9192 & 0.9422 & 0.9468 &
0.9388 & 0.9104 & 0.8406 & 0.6308 & 0.1468 & 0.0038 & 0.9500 & 0.9826\\\hline
$50$ & \multicolumn{1}{r}{$5$} & 0.1614 & 0.6454 & 0.8818 & 0.9542 & 0.9768 &
0.9836 & 0.9878 & 0.9884 & 0.9852 & 0.9788 & 0.9464 & 0.9868 & 0.9958\\
& \multicolumn{1}{r}{$4$} & 0.0384 & 0.3450 & 0.6686 & 0.8206 & 0.8916 &
0.9240 & 0.9396 & 0.9436 & 0.9306 & 0.9026 & 0.8100 & 0.9352 & 0.9728\\
& \multicolumn{1}{r}{$3$} & 0.0088 & 0.0906 & 0.2970 & 0.4702 & 0.6018 &
0.6738 & 0.7164 & 0.7310 & 0.6940 & 0.6340 & 0.4612 & 0.7044 & 0.8236\\
& \multicolumn{1}{r}{$2$} & 0.0034 & 0.0062 & 0.0322 & 0.0860 & 0.1486 &
0.1890 & 0.2256 & 0.2374 & 0.2104 & 0.1766 & 0.0978 & 0.2142 & 0.3206\\
& \multicolumn{1}{r}{$1/2$} & 0.1178 & 0.1728 & 0.2200 & 0.2314 & 0.2268 &
0.1966 & 0.1526 & 0.0918 & 0.0322 & 0.0064 & 0.0030 & 0.2220 & 0.3350\\
& \multicolumn{1}{r}{$1/3$} & 0.5138 & 0.6234 & 0.6906 & 0.7130 & 0.7126 &
0.6722 & 0.6168 & 0.5074 & 0.3090 & 0.0996 & 0.0028 & 0.7110 & 0.8298\\
& \multicolumn{1}{r}{$1/4$} & 0.8424 & 0.8992 & 0.9262 & 0.9336 & 0.9338 &
0.9210 & 0.8978 & 0.8428 & 0.6838 & 0.3744 & 0.0198 & 0.9328 & 0.9698\\
& \multicolumn{1}{r}{$1/5$} & 0.9614 & 0.9796 & 0.9858 & 0.9868 & 0.9868 &
0.9844 & 0.9800 & 0.9632 & 0.8948 & 0.6782 & 0.0958 & 0.9868 & 0.9956\\\hline
$100$ & \multicolumn{1}{r}{$5$} & 0.9958 & 0.9998 & 0.9998 & 1.0000 & 1.0000 &
1.0000 & 1.0000 & 1.0000 & 1.0000 & 1.0000 & 1.0000 & 1.0000 & 1\\
& \multicolumn{1}{r}{$4$} & 0.9502 & 0.9896 & 0.9894 & 0.9986 & 0.9996 &
0.9996 & 0.9996 & 0.9996 & 0.9996 & 0.9996 & 0.9996 & 0.9996 & 1\\
& \multicolumn{1}{r}{$3$} & 0.5996 & 0.8232 & 0.8210 & 0.9586 & 0.9746 &
0.9790 & 0.9832 & 0.9850 & 0.9850 & 0.9800 & 0.9736 & 0.9844 & 0.9946\\
& \multicolumn{1}{r}{$2$} & 0.0460 & 0.1612 & 0.1562 & 0.4158 & 0.5118 &
0.5564 & 0.5926 & 0.6080 & 0.6094 & 0.5716 & 0.5062 & 0.6066 & 0.7594\\
& \multicolumn{1}{r}{$1/2$} & 0.4418 & 0.5162 & 0.4222 & 0.5822 & 0.5906 &
0.5572 & 0.5120 & 0.4420 & 0.3546 & 0.2232 & 0.0900 & 0.6010 & 0.7540\\
& \multicolumn{1}{r}{$1/3$} & 0.9654 & 0.9780 & 0.9632 & 0.9860 & 0.9870 &
0.9830 & 0.9784 & 0.9678 & 0.9474 & 0.8830 & 0.7154 & 0.9874 & 0.9970\\
& \multicolumn{1}{r}{$1/4$} & 0.9996 & 0.9996 & 0.9996 & 0.9996 & 0.9998 &
0.9996 & 0.9996 & 0.9994 & 0.9988 & 0.9962 & 0.9760 & 0.9998 & 1\\
& \multicolumn{1}{r}{$1/5$} & 1.0000 & 1.0000 & 1.0000 & 1.0000 & 1.0000 &
1.0000 & 1.0000 & 1.0000 & 1.0000 & 0.9998 & 0.9980 & 1.0000 & 1\\\hline
$200$ & \multicolumn{1}{r}{$5$} & 1 & 1 & 1 & 1 & 1 & 1 & 1 & 1 & 1 & 1 & 1 &
1 & 1\\
& \multicolumn{1}{r}{$4$} & 1 & 1 & 1 & 1 & 1 & 1 & 1 & 1 & 1 & 1 & 1 & 1 &
1\\
& \multicolumn{1}{r}{$3$} & 0.9996 & 0.9998 & 1.0000 & 1.0000 & 1.0000 &
1.0000 & 1.0000 & 1.0000 & 1.0000 & 1.0000 & 1.0000 & 1.0000 & 1\\
& \multicolumn{1}{r}{$2$} & 0.7552 & 0.8358 & 0.8844 & 0.9026 & 0.9186 &
0.9336 & 0.9378 & 0.9356 & 0.9374 & 0.9350 & 0.9244 & 0.9346 & 0.9742\\
& \multicolumn{1}{r}{$1/2$} & 0.9368 & 0.9400 & 0.9434 & 0.9386 & 0.9312 &
0.9280 & 0.9140 & 0.8818 & 0.8508 & 0.7936 & 0.6930 & 0.9308 & 0.9744\\
& \multicolumn{1}{r}{$1/3$} & 1 & 1 & 1 & 1 & 1 & 1 & 1 & 1 & 1 & 1 & 1 & 1 &
1\\
& \multicolumn{1}{r}{$1/4$} & 1.0000 & 1.0000 & 1.0000 & 1.0000 & 1.0000 &
1.0000 & 1.0000 & 1.0000 & 0.9998 & 0.9986 & 0.9950 & 1.0000 & 1\\
& \multicolumn{1}{r}{$1/5$} & 1 & 1 & 1 & 1 & 1 & 1 & 1 & 1 & 1 & 1 & 1 & 1 &
1\\\hline
\end{tabular}
}
\end{center}
\end{table}

From Tables 3-8, when the simulated critical values are used, we reach the
following conclusions:\newline a) The power increases with $\tau$ ($0<\tau\leq0.5$) and is optimal
when $\tau=0.5$.\newline b) For $\tau=0.2$ the power of the $\widetilde
{LRT}^{(K)}$ is almost always greater than the power of $S^{(K)}$. Moreover,
when $\rho>1$ the performance of $\widetilde{LRT}^{(K)}$ is almost in all
cases the best, while when $\rho<1$ then there is a value of $\lambda$ for
which the performance of $^{0.05}T_{\phi_{\lambda}}^{(K)}$ test statistic is
better. To be more specific it is suggested to use $\lambda=-0.8,-0.9,-1$. At
this point we have to point that when $\rho>1$ it is suggested from Tables 3
and 6 to use $\lambda=-0.1,-0.2,-0.3$. \newline c) The power is increasing
rapidly as the sample size increases.\newline d) For $\tau=0.3,0.5$ the
$S^{(K)}$ has the best performance. Among the test statistics $^{0.05}%
T_{\phi_{\lambda}}^{(K)}$ a similar conclusion is reached with that in
b).\newline e) The results also indicate that when $\tau<0.5$ the test based
on $\widetilde{LRT}^{(K)}$ as well as $S^{(K)}$ performs less good for
$\rho<1$ than in the opposite case of $1/\rho$ (see Haccou et al. (1983, 1985)
for a similar conclusion). However this property does not holds for
$^{0.05}T_{\phi_{\lambda}}^{(K)}$. It is noted that initially there are values
for which the opposite holds and there is a point where the behavior changes.\newline f) When $\tau=0.5$ the test
based on $\widetilde{LRT}^{(K)}$ as well as $S^{(K)}$ performs similar for
$\rho$ and $1/\rho$, with $\rho<1$.\newline

Finally, we have to note some differences and similarities of the results obtained in this paper with the respective ones obtained in Batsidis \textit{et al.} (2011).

Based on Table 2 of this paper, we have concluded that in almost all cases presented here the test statistics are accurate, when the simulated critical values were used. Contrary to this, when the asymptotic critical values of the test statistics were used, it was concluded in Batsidis \textit{et al.} (2011) that for each significance level and sample size
$K\geq40$ there are values of the parameter $\lambda$ for which the power
divergence test statistic is accurate. Moreover, the tests statistics
$\widetilde{LRT}^{(K)}$ and $S^{(K)}$ are in almost all cases not accurate.

Based on the results of this paper and the results of Section 3 of Batsidis \textit{et al.} (2011) related to the power of the test the following common conclusions holds: i) the power is increasing
rapidly as the sample size increases, ii) the power increases with $\tau$ ($0<\tau\leq0.5$) and is optimal
when $\tau=0.5$, and iii) when $\tau<0.5$ the test based
on $\widetilde{LRT}^{(K)}$ as well as $S^{(K)}$ performs less good for
$\rho<1$ than in the opposite case of $1/\rho$ (see Haccou et al. (1983, 1985)
for a similar conclusion). However this property does not holds for
$^{0.05}T_{\phi_{\lambda}}^{(K)}$. It is noted that initially there are values
for which the opposite holds and there is a point where the behavior changes.

In Batsidis \textit{et al.} (2011) it was concluded that there is always a value of $\lambda$ for which $^{0.05}T_{\phi_{\lambda}}^{(K)}$
has greater power than $\widetilde{LRT}^{(K)}$ as well as $S^{(K)}$. This conclusion is not valid when the simulated critical values were used (see conclusion b) above).

\bigskip
\textbf{Acknowledgements. } This work was partially supported by Grant MTM 2009-10072.


\begin{thebibliography}{99}                                                                                               %
\bibitem {Batsidis1}Batsidis, A., Mart\'{\i}n, N., Pardo, L. and Zografos, K.
(2011). Change point analysis of an exponential model based on Phi-divergence test-statistics. Submitted.

\bibitem{car} Cardoso de Oliveira, I.R. and Ferreira, D.F. (2010). Multivariate extension of chi-squared univariate normality test.
\textit{Journal of Statistical Computation and Simulation}, \textbf{80}, 513-526.

\bibitem {hacou2}Haccou, P. , Meelis, E. and Van der Geer, S. (1985). On the
Likelihood Ratio Test for a Change Point in a Sequence of Independent
Exponentially Distributed Random Variables., \emph{Report MS-R8507, Centre for
Mathematics and Computer Science, }Amsterdam, The Netherlands.

\bibitem {hacou1}Haccou, P., Meelis, E. and Van der Geer, S. (1988). The
Likelihood Ratio Test for the Change point Problem for Exponentially
distributed Random Variables.\textit{ }\emph{Stochastic Processes and their
Applications, }\textbf{27, }121-139.

\bibitem {HorvathSerb}Horv\'{a}th, L. and Serbinowska, M. (1995). Testing for
Changes in Multinomial Observations: the Lindisfarne Scribes problem.
\emph{Scandinavian Journal of Statistics}, \textbf{22}, 371--384.

\bibitem {Pardo}Pardo, L. (2006). \emph{Statistical inference based on
divergence measures}. Chapman \& Hall/CRC, Boca Raton.

\bibitem {rom}Romeu, Jorge Luis and Ozturk, Aydin (1993). A comparative study
of goodness-of-fit tests for multivariate normality. \textit{J. Multivariate
Anal.}, \textbf{46 }, 309--334.

\bibitem {srihui}Srivastava, M. S. and Hui, T. K. (1987). On assessing
multivariate normality based on Shapiro-Wilk $W$ statistic. \textit{Statist.
Probab. Lett.}, \textbf{5},15--18.

\bibitem {VostricovaJu}Vostrikova, L. Ju. (1981). Detecting disorder in
multidimensional random processes. \emph{Soviet Math Dokl}. \textbf{24}, 55--59.

\bibitem{wor}  Worsley, K. J. (1986) Confidence regions and test for a change-point in a sequence of exponential family random variables. \emph{Biometrika}, \textbf{73}, 91�104.
\end{thebibliography}
\end{document}